\documentclass{article}[11pt]
\usepackage{amsmath,amssymb,amsthm,a4wide}
\usepackage{amssymb}
\input xypic
\usepackage[all,tips]{xy}
\usepackage{graphicx}               

\usepackage[hidelinks]{hyperref}

\usepackage{tikz}
\usetikzlibrary{decorations.markings}
\usetikzlibrary{arrows}
\usetikzlibrary{arrows.meta}

\tikzset{middleprightarrow/.style={decoration={markings,mark= at position 0.5 with {\arrow[scale=1.5,>=stealth]{>}},},postaction={decorate}}}
\tikzset{middlepleftarrow/.style={decoration={markings,mark= at position 0.5 with {\arrow[scale=1.5,>=stealth reversed]{>}},},postaction={decorate}}}

\tikzset{rightprightarrow/.style={decoration={markings,mark= at position 0.75 with {\arrow[scale=1.5,>=stealth]{>}},},postaction={decorate}}}
\tikzset{rightpleftarrow/.style={decoration={markings,mark= at position 0.75 with {\arrow[scale=1.5,>=stealth reversed]{>}},},postaction={decorate}}}

\tikzset{leftprightarrow/.style={decoration={markings,mark= at position 0.25 with {\arrow[scale=1.5,>=stealth]{>}},},postaction={decorate}}}
\tikzset{leftpleftarrow/.style={decoration={markings,mark= at position 0.25 with {\arrow[scale=1.5,>=stealth reversed]{>}},},postaction={decorate}}}

\newtheorem{thm}{Theorem}[section]
\newtheorem{lem}[thm]{Lemma}
\newtheorem{prop}[thm]{Proposition}
\newtheorem{cor}[thm]{Corollary}
\newtheorem{conj}[thm]{Conjecture}

\theoremstyle{definition}
\newtheorem{defin}[thm]{Definition}

\theoremstyle{remark}
\newtheorem{remark}[thm]{Remark}
\newtheorem{example}[thm]{Example}

\newcommand{\bth}{\begin{thm}}
\renewcommand{\eth}{\end{thm}}
\newcommand{\bpr}{\begin{prop}}
\newcommand{\epr}{\end{prop}}
\newcommand{\ble}{\begin{lem}}
\newcommand{\ele}{\end{lem}}
\newcommand{\bco}{\begin{cor}}
\newcommand{\eco}{\end{cor}}
\newcommand{\bde}{\begin{defin}}
\newcommand{\ede}{\end{defin}}
\newcommand{\bex}{\begin{example}}
\newcommand{\eex}{\end{example}}
\newcommand{\bre}{\begin{remark}}
\newcommand{\ere}{\end{remark}}
\newcommand{\bcj}{\begin{conj}}
\newcommand{\ecj}{\end{conj}}

\newcommand{\beq}{\begin{equation}}
\newcommand{\eeq}{\end{equation}}

\newcommand{\ot}{{\otimes}}

\newcommand{\lb}{\label}
\newcommand{\nl}{\newline}
\newcommand{\bpf}{\begin{proof}}
\newcommand{\epf}{\end{proof}}
\newcommand{\da}{{\text -}}

\newcommand{\E}{{\cal E}}

\newcommand{\C}{{\cal C}}

\newcommand{\D}{{\cal D}}

\newcommand{\caL}{{\cal L}}
\newcommand{\cS}{{\cal S}}
\renewcommand{\S}{{\cal S}}

\newcommand{\M}{{\cal M}}

\newcommand{\Rep}{{\cal R}{\it ep}}
\newcommand{\Set}{{\cal S}{\it et}}
\newcommand{\Vect}{{\cal V}{\it ect}}

\newcommand{\Del}{{\bf \Delta}}
\newcommand{\gl}{{\mathfrak g\mathfrak l}}
\newcommand{\g}{{\mathfrak g}}
\newcommand{\z}{{\mathfrak z}}
\newcommand{\SW}{{\cal SW}}

\newcommand\void[1]{}

\begin{document}
\author{Alexei Davydov $^a)$ and Mohamed Elbehiry $^b)$}
\title{Deformation cohomology of Schur-Weyl categories. \\ Free symmetric categories}
\maketitle
\date{}
\maketitle
\begin{center}
$a)$ Department of Mathematics, Ohio University, Athens, OH 45701, USA\\
$b)$ Department of Mathematics, Northeastern University, Boston, MA 02115, USA\\
\end{center}
\date{}
\begin{abstract}
The deformation cohomology of a tensor category controls deformations of its monoidal structure. 
Here we describe the deformation cohomology of tensor categories generated by one object (the so-called Schur-Weyl categories).
Using this description we compute the deformation cohomology of free symmetric tensor categories generated by one object with an algebra of endomorphism free of zero-divisors. 
We compare the answers with the exterior invariants of the general linear Lie algebra. 
\end{abstract}
\tableofcontents

\section{Introduction}

Possible monoidal structures on a given tensor category naturally form an object of an algebro-geometric nature (the moduli space).
The tangent space to the moduli space of tensor structures is computed by the third cohomology of a certain complex, the deformation complex of the tensor category \cite{da, ye}\footnote{the deformation cohomology of a tensor category is sometimes called its Davydov-Yetter cohomology.}. 

The aim of this paper is to compute the deformation cohomology of free symmetric tensor categories generated by one object. 
Free symmetric tensor categories generated by one object have a property that there are non nontrivial morphisms between different tensor powers of the generator. Such tensor categories were dubbed Schur-Weyl categories in \cite{dm}.
We explore the internal organisation of the deformation complex of a Schur-Weyl category. 
This complex comes equipped with a natural decreasing filtration. Component of the associated graded complex (called horizontal complexes) have a uniform structure. Horizontal complexes can be presented as cochain complexes associated with diagrams of vector spaces in the form of a higher-dimensional cube (called cubic diagrams here). 
Cubic diagrams coming from free symmetric tensor categories are diagrams of invariants ot a symmetric group representation.
Using the relation with the cohomology of simplicial cubes we compute the cohomology of cubic diagrams of invariants. 
This allows us to describe the deformation cohomology of free symmetric tensor categories. 
Let $\S$ be the free symmetric tensor category generated by one object whose endomorphisms are scalar. 
We show (assuming that the ground field is of odd characteristic) that the deformation cohomology $H^*(\cS)$ is the exterior algebra 
$\Lambda(e_1,e_3,e_5,...)$ on odd degree generators $deg(e_{2i-1}) = 2i-1$ (theorem \ref{cfs}). 
Let $\S(A)$ be the free symmetric tensor category generated by one object whose endomorphism algebra $A$ is a domain. 
Then (assuming that $A$ is a commutative algebra, which is more than one dimensional)
the deformation cohomology of the free symmetric category $\cS(A)$ is the exterior algebra of the first cohomology $H^1(\cS(A)) = A$, i.e.
$H^*(\cS(A))\ \simeq\ \Lambda^*(A)$ (theorem \ref{cfsa}). 
We also use our methods to compute the deformation cohomology of the degenerate affine Hecke category $\caL$, i.e. the Schur-Weyl category corresponding to the collection of the degenerate affine Hecke algebras. We show that $H^*(\caL)$ is the exterior algebra $\Lambda^*(k[x])$  of the algebra of polynomials (theorem \ref{dvda}). 

The free symmetric tensor category $\cS(A)$ can be thought of as the limiting case of the representation category $\Rep(\gl(V)\ot A)$ of the general linear Lie algebra $\gl(V)\ot A$, when the dimension of the vector space $V$ goes to infinity. 
The deformation cohomology of the representation category $\Rep(\g)$ can be identified with the adjoint $\g$-invariants of the exterior algebra $\Lambda^*(\g)^\g$ (theorem \ref{dcrc}, see also \cite{da}). Here we assume that the characteristic of the ground field is zero.
Classical invariant theory says that $\Lambda^*(\gl(V))^{\gl(V)}$ is he exterior algebra $\Lambda(x_1,x_3,...,x_{2d-1})$ with generators of degree $deg(x_{2i-1}) = 2i-1$ and with $d=dim(V)$ (see e.g. \cite{ko,it}).
We use the Schur-Weyl duality functor $SW:\S\to \Rep(\gl(V))$, sending the generator to the vector representation $V$, to relate the deformation cohomology of $\S$ and of $\Rep(\gl(V))$. 
Note that the functoriality property of the deformation cohomology is not straightforward and is similar to the functoriality of the centre, or the Hochschild cohomology (see \cite{dkr}): a tensor functor $F:\C\to\D$ gives rise to a cospan of homomorphisms of graded algebras
$$\xymatrix{ & H^*(F) \\ H^*(\C) \ar[ur] && H^*(\D) \ar[ul]}$$
where $H^*(F)$ is the deformation cohomology of tensor functor $F$. 
We show that the deformation cohomology of the Schur-Weyl functor $SW$ is the exterior algebra 
$H^*(SW) = \Lambda(e_1,e_3,...,e_{2d-1})$,  
the homomorphism $H^*(\S)\ \to\ H^*(SW)$ is the quotioning by the ideal generated by $e_{s}, s>2d-1$, and 
the homomorphism $H^*(\Rep(\gl(V)))\ \to\ H^*(SW)$ is an isomorphism sending $x_m$ to $((m-1)!)^{-1}e_{m}$ 
(theorem \ref{cswf}).
This in particular give a precise formulation to the intriguing connection between the combinatorics of partitions (giving the answer for $H^*(\S)$) and the exterior invariants of $\gl(V)$ observed by Kostant in \cite{ko}. 
We also relate the deformation cohomology of $\S(A)$ and of $\Rep(\gl(V\ot A))$ under the assumption that $A$ be a commutative algebra with $dim(A)>1$ and such that $A^{\ot n}$ has no zero-divisors for any $n$. 
Then the deformation cohomology of the Schur-Weyl functor $SW:\S(A)\to \Rep(\gl(V\ot A))$ (sending the generator to the vector representation $V\ot A$) is the exterior algebra $H^*(SW) = \Lambda^*(A)$ and the homomorphisms 
$$H^*(\S(A))\ \to\ H^*(SW)\ \leftarrow\ H^*(\Rep(\gl(V\ot A)))$$ 
are isomorphisms (theorem \ref{swda}).

The one-dimensional cohomology $H^3(\S)$ suggests that the moduli space of tensor structures of $\S$ is (locally) one-dimensional.
This was shown to be true in \cite{da}. The detailed analysis of this moduli space and its relation to the one-parameter family of Hecke categories from \cite{dm} will be the subject of a subsequent publication. 

Throughout $k$ will be the ground field. 
We denote by $\Vect$ the category of vector spaces over $k$. 
By a tensor category over $k$ we mean a monoidal category enriched over $\Vect$. 
A tensor functor between tensor categories is a $\Vect$-enriched monoidal functor.

\subsection*{Acknowledgment}

This paper had a rather difficult early life. Started in 2016, when the second author was a master student at Ohio University it was largely complete by the time the second author moved to Northeastern University. Mostly due to the ineffectiveness of the first author in managing his ever increasing load it took three years to do the final polishing. 
The first author would like to thank Max Planck Institute for Mathematics (Bonn, Germany) for hospitality during the Summer of 2019, crucial for the completion of this paper. 

\section{Deformation cohomology of tensor categories}

Let $\C,\D$ be  tensor categories, which we assume to be strict for simplicity of exposition.
For a tensor functor $F:\C\to\D$ define its $n$-th power by
$$F^{\ot n}:\C\times ...\times\C\to \D,\qquad F^{\ot n}(X_1,...,X_n) = F(X_1 \ot  ...\ot X_{n})$$
For $n=0$ denote
$$F^{\otimes 0}:\Vect \to\D,\qquad F^{\otimes 0}(V) = V\otimes I $$
where $I\in\D$ is the unit object.

\subsection{Deformation complex of a tensor functor}

Denote by $E^*(F) = End(F^{\otimes *})$ the collection of endomorphisms algebras of tensor powers of a monoidal functor $F$.
An element of $E^n(F)$ is a collection of endomorphisms $a_{X_{1},...,X_{n}}\in End_\D(F(X_1 \ot  ...\ot X_{n}))$ for $X_{1},...,X_{n}\in\C$ natural in the following sense: 
the diagram 
\beq\lb{nat}\xymatrix{
F(X_1 \ot  ...\ot X_{n}) \ar[rr]^{a_{X_{1},...,X_{n}}} \ar[d]_{F(f_1\ot...\ot f_n)} && F(X_1 \ot  ...\ot X_{n}) \ar[d]^{F(f_1\ot...\ot f_n)} \\
F(Y_1 \ot  ...\ot Y_{n}) \ar[rr]^{a_{Y_{1},...,Y_{n}}} && F(Y_1 \ot  ...\ot Y_{n}) 
}\eeq
commutes for any $f_i\in \C(X_i,Y_i)$.
\nl
Here we recall (following \cite{da}) how to equip the collection of algebras $E^{*}(F)$ with the structure of cosimplicial algebra.
\nl
More precisely the image of the coface map
$${\partial}_{i} : End(F^{\otimes n}) \rightarrow End(F^{\otimes n+1}) \qquad i = 0,...,n+1$$
of an endomorphism $a \in End(F^{\otimes n})$ has the following specialisation on objects $X_{1},...,X_{n+1}\in\C$:
$${\partial}_{i}(a )_{X_{1},...,X_{n+1}} = \left\{
\begin{array}{lcc}
\phi(1_{F(X_{1})} \otimes a_{X_{2},...,X_{n+1}})\phi^{-1} & , & i = 0 \\
a_{X_{1},...,X_{i} \otimes X_{i+1},...,X_{n+1}} & , & 1 \leq i \leq n \\
\psi(a_{X_{1},...,X_{n}} \otimes 1_{F(X_{n+1})})\psi^{-1} & , & i = n+1
\end{array}
\right.$$
Here $\phi$ is the tensor structure constraint  
$F(X_1) \ot F(X_2\ot ...\ot X_{n+1})\ \to\ F(X_1 \ot  ...\ot X_{n+1})\ ,$
and $\psi$ is the tensor structure constraints 
$F(X_1\ot  ...\ot X_{n+1})\ \to\  F(X_1 \ot  ...\ot X_{n})\ot F(X_{n+1})\ .$
\nl
The specialisation of the image of the codegeneration map
$${\sigma}_{i} : End(F^{\otimes n}) \rightarrow End(F^{\otimes n+1}) \qquad i = 0,...,n-1$$
is
$${\sigma}_{i}(a)_{X_{1},...,X_{n-1}} = {a}_{X_{1},...,X_{i},I,X_{i+1},...,X_{n-1}}.$$
The zero component of this complex is the  endomorphism algebra $End_\D(I)$ of the unit object $I$ of the category $\D$, which can be regarded as the endomorphism algebra of the functor $F^{\otimes 0}$.
\nl
The coface maps
$${\partial}_{i}:End_\D(I) \rightarrow End(F),\qquad i=0,1$$
have the form
$${\partial}_{0}(a)_X = \rho_{F(X)} (a\otimes 1_{F(X)}){\rho}_{F(X)}^{-1}, \qquad {\partial}_{1}(a)_X = \lambda_{F(X)} (1_{F(X)}\otimes a){\lambda_{F(X)}}^{-1};$$
here $\rho_{F(X)}:I\ot {F(X)}\to {F(X)}$ and $\lambda_{F(X)}:{F(X)}\ot I\to {F(X)}$ are the structural isomorphisms of the unit object $I$.

It is straightforward to verify the following.
\bpr
The maps ${\sigma}_i$ and ${\partial}_j$ make $E^*(F)$ a cosimplicial complex.
\epr

The cohomology of the corresponding cochain differential 
\beq\lb{chdi}
\partial: E^n(F)\to E^{n+1}(F),\qquad \partial = \sum_{i=0}^{n+1}(-1)^i{\partial}_{n+1}^{i}
\eeq
was called the {\em  tangent cohomology} of $F$ in \cite{da}. Here we call it the {\em deformation cohomology} of $F$. 

\bex
The space of 1-cocycles $Z^1(F)$ coincides with the space
$$Der(F) = \{a\in End(F)|\  F_{X,Y}^{-1}a_{X\ot Y}F_{X,Y} = 1_{F(X)}\ot a_Y+ a_X\ot 1_{F(Y)},\quad X,Y\in \C\}$$
of {\em derivations} (or {\em primitive} endomorphisms) of $F$.
\nl
The subspace of 1-coboundaries $B^1(F)\subset Z^1(F)$ corresponds to the subspace $Der_{inn}(F)$ of {\em inner derivations} of $F$.
The first cohomology $H^1(F)$ is the space $OutDer(F) = Der(F)/Der_{inn}(F)$ of {\em outer derivations} of $F$.
\eex

The deformation complex $E^*(F)$ is equipped with the {\em $\cup$-product} 
$$\cup:E^m(F)\ot E^n(F)\ \to\ E^{m+n}(F)$$
$$(a\cup b)_{X_1,...,X_m,X_{m+1},...,X_{m+n}} = \phi(a_{X_1,...,X_m}\ot b_{X_{m+1},...,X_{m+n}})\phi^{-1}\ ,\qquad a\in E^m(F),\  b\in E^n(F) \ .$$
Here $\phi= F_{X_1\ot ... \ot X_{m},X_{m+1} \ot ...\ot X_{m+n}}$ is the 
coherence isomorphism
$$F(X_1 \ot  ... \ot X_{m})\ot F(X_{m+1} \ot ... \ot X_{m+n})\to  F(X_1 \ot  ... \ot X_{m+n})\ .$$
The $\cup$-product induces an associative multiplication on the cohomology
\beq\lb{cup}
\cup: H^m(F)\ot H^n(F)\ \to\ H^{m+n}(F)\ .
\eeq
Methods similar to those from \cite{ge} (see also \cite{bd}) show that the cup-product is graded commutative
$$b\cup a = (-1)^{|a||b|}a\cup b\qquad a\in H^m(F), \ b\in H^n(F) \ .$$
Indeed, the commutativity homotopy for the $\cup$-product can be chosen as
$$a* b =  \sum_{i=1}^m(-1)^{(n-1)i}a*_ib\ ,$$
where 
$$(a*_ib)_{X_1,...,X_{m+n-1}} = $$
$$= (1_{F(X_1\ot...\ot X_{i-1})}\ot b_{X_{i},...,X_{i+m-1}}\ot 1_{F(X_{i+m}\ot...\ot X_{m+n})})a_{X_1,...,X_{i-1},X_{i}\ot...\ot X_{i+m-1},X_{i+m},...,X_{m+n-1}}\ .$$
\\

Let now $F:\C\to\D$ and $G:\D\to\E$ be tensor functors.
We have two collections of algebra homomorphisms
\beq\lb{cma}E^n(G)\ \to\ E^{n}(G\circ F)\ ,\qquad E^n(F)\ \to\ E^{n}(G\circ F)\ .\eeq
The first sends $a\in E^n(G)$ into 
$$\xymatrix{G(F(X_1 \ot  ...\ot X_{n})) \ar[d]_{G(\psi)^{-1}} & && G(F(X_1 \ot  ...\ot X_{n}))  \\
G(F(X_1) \ot  ...\ot F(X_{n})) \ar[rrr]^{a_{F(X_1),  ..., F(X_{n})}} &&& G(F(X_1) \ot  ...\ot F(X_{n}))  \ar[u]_{G(\psi)} }$$
where $\psi:F(X_1) \ot  ...\ot F(X_{n})\to F(X_1 \ot  ...\ot X_{n})$ is the tensor constraint of $F$,
while the second sends $b\in E^n(F)$ into 
$$\xymatrix{G(F(X_1 \ot  ...\ot X_{n}))  \ar[rrr]^{G(b_{X_1,  ..., X_{n}})} &&&  G(F(X_1 \ot  ...\ot X_{n}))}\ .$$
It is straightforward to see that the homomorphisms \eqref{cma} are cosimplicial maps and give rise to homomorphisms of deformation cohomology
\beq\lb{hma}H^n(G)\ \to\ H^{n}(G\circ F)\ ,\qquad H^n(F)\ \to\ H^{n}(G\circ F)\ .\eeq
\\

We are mostly interested in the case $F=Id_\C$. We denote $E^*(\C)=E^*(Id_\C), \ H^*(\C)=H^*(Id_\C)$ and call them the {\em deformation complex} and the {\em deformation cohomology} of the tensor category $\C$. 
Using the construction \eqref{hma} a tensor functor $F:\C\to\D$ gives rise to a pair of homomorphisms
\beq\lb{hmc}H^n(\C)\ \to\ H^{n}(F)\ ,\qquad H^n(\D)\ \to\ H^{n}(F)\ .\eeq

\ble\lb{cef}
Let $F:\C\to\D$ be a fully faithful tensor functor. Then the homomorphism $H^*(\C)\ \to\ H^*(F)$ is an isomorphism.
\nl
Let $F:\C\to\D$ be a full tensor functor essentially surjective on objects. Then the homomorphism $H^*(\D)\ \to\ H^*(F)$ is an isomorphism.
\ele
\bpf
Let $F:\C\to\D$ be a fully faithful tensor functor. To show bijectivity of $E^n(\C)\to E^n(F)$ note that for any $X_1,...,X_n\in\C$ the effect on morphisms $End_\C(X_1\ot...\ot X_n)\to End_\D(F(X_1\ot...\ot X_n))$ is an isomorphism. Moreover the naturality conditions \eqref{nat} for $Id_\C$ and $F$ are the same.
\nl
Let now $F:\C\to\D$ be a full tensor functor essentially surjective on objects. To show bijectivity of $E^n(\D)\to E^n(F)$ note that for any $Y_1,...,Y_n\in\C$ there are $X_1,...,X_n\in\C$ such that $F(X_i)\simeq Y_i$. Thus we have isomorphisms 
$$End_\D(Y_1\ot...\ot Y_n)\to End_\D(F(X_1)\ot...\ot F(X_n))\to End_\D(F(X_1\ot...\ot X_n))\ .$$
Using the surjections $\C(X_i,X_i')\to \D(F(X_i),F(X_i'))\simeq \D(Y_i,Y_i')$ we can see that the naturality conditions \eqref{nat} for $Id_\D$ and $F$ are identified by these isomorphisms.
\epf
Together with \eqref{hmc} the above lemma assigns a homomorphism $H^*(\D)\to H^*(\C)$ to a fully faithful tensor functor $F:\C\to\D$ and a homomorphism $H^*(\C)\to H^*(\D)$ to a full tensor functor $F:\C\to\D$ essentially surjective on objects.

\section{Schur-Weyl categories and their deformation cohomology}

\subsection{Schur-Weyl categories}

Here we recall the definitions from \cite{dm} and set our notations.
\nl
A {\em multiplicative sequence of algebras} is 
a collection of associative unital algebras $A_*=\{A_n\ |\ n\geqslant 0\}$ (with $A_0=k$) 
equipped with a collection of (unital) algebra homomorphisms
$$\mu_{m,n}:A_m\otimes A_n\to A_{m+n},\qquad m,n\geqslant 0,$$
satisfying the following associativity axiom: for any $l,m,n\geqslant 0$ the diagram
\beq\lb{asam}
\xymatrix{ A_l\otimes A_m\otimes A_n
\ar[rr]^{\mu_{l,m}\ot I} \ar[d]_{I\ot\mu_{m,n}} &&
A_{l+m}\otimes A_n \ar[d]^{\mu_{l+m,n}} \\
A_l\otimes A_{m+n} \ar[rr]^{\mu_{l,m+n}} && A_{l+m+n} .}
\eeq
commutes. 
\nl
We denote by $\mu_{m_1,...,m_n}:A_{m_1}\ot...\ot A_{m_n}\ \to\ A_{m_1+...+m_n}$ the unique composition of homomorphisms $\mu_{m,n}$.

A model example of a multiplicative sequence of algebras is provided by the following construction.
\bex
Let $\C$ be a (strict) tensor category
such that $End_\C(I)=k$,
where $I$ denotes the unit object of $\C$.
Given an object $X$ of
$\C$, the sequence $A_*$ with $A_n=End_\C(X^{\otimes n})$
is multiplicative
with respect to the homomorphisms $\mu_{m,n}$ given by the tensor
product on morphisms
$$
End_\C(X^{\otimes m})\otimes End_\C(X^{\otimes n}) \to End_\C(X^{\otimes m+n}).
$$
\eex
Moreover, any multiplicative sequence can be
obtained in this way. Indeed,
starting with a multiplicative sequence $A_*$, define its {\em Schur-Weyl category}
$\SW(A_*)$ with objects $[n]$ parameterized by natural numbers,
with no morphisms between different objects and with the endomorphism
algebras $End_{\SW(A_*)}([n])=A_n$. Define tensor product
on the objects of $\SW(A_*)$ by $[m]\otimes [n]=[m+n]$.
The multiplicative structure of the sequence $A_*$ yields the tensor
product on morphisms:
$$
\xymatrix{End_{\SW(A_*)}([m])\otimes End_{\SW(A_*)}([n]) = A_m\ot A_n \ar[r]^(.58){\mu_{m,n}} & A_{m+n}= End_{\SW(A_*)}([m+n])}.
$$
Note that the Schur-Weyl category $\SW(A_*)$ is a strict (and skeletal) tensor category.

\bre\lb{nstr}
Let $\C$ be a non-strict monoidal linear category with the associativity constraint $\alpha$. For $X\in\C$ define $X^{\ot n}$ inductively by $X^{\ot n} = X\ot X^{\ot n-1}$ with $X^{\ot 0} = I$. 
Define isomorphisms $\alpha_{m,n}:X^{\ot m}\ot X^{\ot n}\to X^{\ot m+n}$ inductively by
$$\xymatrix{X^{\ot m}\ot X^{\ot n}\ \ar[rrrr]^{\alpha_{m,n}} \ar@{=}[d] &&&& {\quad X^{\ot m+n}} \\
(X\ot X^{\ot m-1})\ot X^{\ot n} \ar[rr]^{\alpha_{X,X^{\ot m-1},X^{\ot n}}} && X\ot (X^{\ot m-1}\ot X^{\ot n}) \ar[rr]^{1\ot \alpha_{m-1,n}} && {\quad X\ot X^{\ot m+n-1}} \ar@{=}[u]}
$$
The homomorphisms $$
\mu_{m,n}:End_\C(X^{\otimes m})\otimes End_\C(X^{\otimes n}) \to End_\C(X^{\otimes m+n}).
$$
can now be defined by
$$\mu_{m,n}(a,b) = \alpha_{m,n}(a\ot b)\alpha_{m,n}^{-1}\ .$$
It follows from Mac Lane coherence that the homomorphisms $\mu_{m,n}$ make $End_\C(X^{\otimes n})$ a multiplicative sequence of algebras.
\ere

\bre
The notation used for this construction in \cite{dm} was $\overline\SW(A_*)$ and $\SW(A_*)$ was used for the $k$-linear abelian envelope of $\overline\SW(A_*)$. Since here $\overline\SW(A_*)$ is the more basic object we simplify its notation. 
\ere

Let $f_*=\{f_n\ |\ n\geqslant 0\}$ be a sequence of
algebra homomorphisms
$f_n:A_n\to B_n$ between the corresponding algebras of two
multiplicative sequences $A_*$ and $B_*$.
We call $f_*$ a {\it homomorphism
of multiplicative sequences} if for any $m,n$ the following diagram
commutes:
\beq\lb{hms}\xymatrix{ A_m\otimes A_n \ar[rr]^{f_m\otimes f_n}
\ar[d]_{\mu_{m,n}} && B_m\otimes B_n \ar[d]^{\mu_{m,n}} \\
A_{m+n} \ar[rr]^{f_{m+n}} && B_{m+n} .}
\eeq
We will say that $f_*$ is an {\it epimorphism\/}, if all homomorphisms $f_n$ are surjective.

A model construction of a homomorphism of multiplicative sequences is provided by the following
\bex
Let $F:\C\to\D$ be a tensor functor between (strict) tensor categories (with $End_\C(I)=End_\D(I)=k$).
Given an object $X\in\C$,  define a sequence of homomorphisms $f_*:A_*\to B_*$ with $A_n=End_\C(X^{\otimes n})$ and $B_n=End_\D(F(X)^{\otimes n})$ by 
$$f_n(a) = \phi_n(F(a))\phi_n^{-1}\ ,$$
where isomorphisms $\phi_n:F(X^{\ot n})\to F(X)^{\ot n}$ are defined inductively by
$$\xymatrix{F(X^{\ot n})\ \ar[rrrr]^{\phi_{n}} \ar@{=}[d] &&&& {\quad F(X)^{\ot n}} \\
F(X\ot X^{\ot n-1}) \ar[rr]^{F_{X,X^{\ot n-1}}} && F(X)\ot F(X^{\ot n-1}) \ar[rr]^{1\ot \phi_{n-1}} && {\quad F(X)\ot F(X)^{\ot n-1}} \ar@{=}[u]}
$$
\eex
Moreover, any homomorphism of multiplicative sequences can be obtained in this way.
Indeed, the construction of Schur-Weyl categories is functorial with respect to homomorphisms of
multiplicative sequences: a homomorphism $f_*:A_*\to B_*$ defines a strict tensor functor $F(f_*):\SW(A_*)\to\SW(B_*)$.

A {\em gaude transformation} between homomorphisms of multiplicative sequences $f_*, g_*:A_*\to B_*$ is a collection invertible elements $c(n)\in B_{n}^\times$ such that $c(n)f_n(a) = g_n(a)c(n)$ for any $a\in A_n$ and such that
$$c(m+n) = \mu_{m,n}(c(m)\ot c(n)),\qquad \forall m,n\ .$$

The following is straightforward.
\bpr\lb{twh}
Strict tensor functors $F:\SW(A_*)\to\SW(B_*)$ between Schur-Weyl categories such that  $F([1])=[1]$ are in 1-1 correspondence with homomorphisms $A_*\to B_*$ of multiplicative sequences.
\nl
Tensor natural transformations between tensor functors corresponding to twisted homomorphisms $f_*, f(*,*)$ and $g_*, g(*,*)$ are gaude transformation between these twisted homomorphisms. 
\epr

\subsection{Deformation cohomology of Schur-Weyl categories}

\bpr
The deformation complex of the Schur-Weyl category $\SW(A_*)$ is
$$E^n(\SW(A_*)) = \bigoplus_{m_1,...,m_n} C_{A_{m_1+...+m_n}}(A_{m_1}\ot...\ot A_{m_n})\ ,$$
where the direct sum is taken over $n$-tuples of positive integers.
\nl
The cosimplicial differentials restricted to $C_{A_{m_1,...,m_n}}(A_{m_1}\ot...\ot A_{m_n})$ are given by the following embeddings:
\nl
the zero differential is the direct sum (over $m_0$) of 
$$C_{A_{m_1+...+m_n}}(A_{m_1}\ot...\ot A_{m_n})\ \to\ C_{A_{m_0+m_1+...+m_n}}(A_{m_0}\ot A_{m_1}\ot...\ot A_{m_n})\ ;$$
for $0<i<n+1$ the differential $\partial_i$ is the direct sum (over $m_i'$ and $m_i'$ such that $m_i'+m_i''=m_i$ ) of
$$C_{A_{m_1+...+m_n}}(A_{m_1}\ot...\ot A_{m_n})\ \to\ C_{A_{m_1+...m_i'+m_i"+...+m_n}}(A_{m_1}\ot...A_{m_i'}\ot A_{m_i''}\ot ...\ot A_{m_n})\ ;$$
the differential $\partial_{n+1}$ is the direct sum (over $m_{n+1}$) of 
$$C_{A_{m_1+...+m_n}}(A_{m_1}\ot...\ot A_{m_n})\ \to\ C_{A_{m_1+...+m_n+m_{n+1}}}(A_{m_1}\ot...\ot A_{m_n}\ot A_{m_{n+1}})\ .$$
\epr
\bpf
By the definition $E^n(\SW(A_*))$ is the space of natural in $X_1,...,X_n\in \SW(A_*)$ collections of endomorphisms $a_{X_1,...,X_n}:X_1\ot...\ot X_n\to X_1\ot...\ot X_n$.
By naturality any such collection is determined by the collection of its values for $X_1=[m_1],...,X_n=[m_n]$. 
Being an endomorphism of $[m_1]\ot...\ot[m_n]=[m_1+...+m_n]$ a value $a_{[m_1],...,[m_n]}$ is an element of $A_{m_1+...+m_n}$. 
Naturality of the collection $a$ is equivalent to the condition that values $a_{[m_1],...,[m_n]}$ belong to the centralisers $C_{A_{m_1+...+m_n}}(A_{m_1}\ot...\ot A_{m_n})$. 

\epf

Denote 
\beq\lb{cen}C(m_1,...,m_n) = C_{A_{m_1+...+m_n}}(A_{m_1}\ot...\ot A_{m_n})\ .\eeq
Below is the picture of the first four layers of the deformation cohomology complex.
$$\xygraph{ !{0;/r10pc/:;/u6pc/::}
*+{C(1)}
(
:@{}[d] *+{C(2)} 
 (
 :@{}[d] *+{C(3)}
  (
  :@{}[d] *+{C(4)}
   (
   :[u(.3)r(1.2)] *+{C(1,3)}="13"
    (
    :[r] *+{C(1,1,2)}="112"
    :[d(.3)r(.8)] *+{C(1,1,1,1)}="1111"
    ,
    :[d(.3)r(.8)] *+{C(1,2,1)}="121"
    :"1111"
    )
   ,
   :[r] *+{C(2,2)}="22"
    (
    :"112"
    ,
    :[d(.3)r(.8)] *+{C(2,1,1)}="211"
    :"1111"
    )
   ,
   :[d(.3)r(.8)] *+{C(3,1)}="31"
    (
    :"211"
    ,
    :"121"
    )
   )
  ,
  :"13"
  ,
  :"31"
  ,
  :[u(.3)r(1.2)] *+{C(1,2)}="12"
   (
   :"112"
   ,
   :"121"
   ,
   :[d(.3)r(.8)] *+{C(1,1,1)}="111"
   :"1111"
   )
  ,
  :[d(.3)r(.8)] *+{C(2,1)}="21"
   (
   :"111"
   ,
   :"121"
   ,
   :"211"
   )
  )
 ,
 :"21"
 ,
 :"12"
 ,
 :[r] *+{C(1,1)}="11"
  (
  :"111"
  ,
  :"112"
  ,
  :"211"
  )
  ,
:"22"
 )
,
:"11"
,
:"21"
,
:"12"
,
:"13"
,
:"31"
)
}$$
Arrows represent components of the differential.
Arrows in a layer (homomorphisms between centralisers $C(m_1,...,m_n)$ preserving the sum $m_1+...+m_n$) correspond to the differentials $\partial_i$ for $0<i<n+1$. 
Arrows between layers are components of the differentials $\partial_0$ and $\partial_{n+1}$.
Note that some of those arrows correspond to two different homomorphisms, e.g. $C(1,1)\to C(1,1,1)$ represents both $\partial_0$ and $\partial_3$.
\\

Define a decreasing filtration on $E^*(\SW(A_*))$ by 
\beq\lb{fil}F^pE^n(\SW(A_*)) = \bigoplus_{m_1+...+m_n> p} C_{A_{m_1+...+m_n}}(A_{m_1}\ot...\ot A_{m_n})\ .\eeq
Clearly $F^pE^n(\SW(A_*))$ is a subcomplex of $E^n(\SW(A_*))$.
We call the associated graded complex $E^n_p=F^pE^n/F^{p+1}E^n$ the {\em horizontal complex of depth} $p$. 

The filtration \eqref{fil} gives rise to the spectral sequence
\beq\lb{ss}E^{p,q}_1 = H^p(E^*_{p+q})\ \Rightarrow \ H^{p+q}(\SW(A_*))\eeq
converging to the deformation cohomology $E^n(\SW(A_*))$.

The first four horizontal complexes are depicted below.
$$\xygraph{ !{0;/r10pc/:;/u6pc/::}
*+{C(1)}
(
:@{}[d] *+{C(2)} 
 (
 :@{}[d] *+{C(3)}
  (
  :@{}[d] *+{C(4)}
   (
   :[u(.3)r(1.2)] *+{C(1,3)}="13"  ^{-}
    (
    :[r] *+{C(1,1,2)}="112"
    :[d(.3)r(.8)] *+{C(1,1,1,1)}="1111"  ^(.4){-}
    ,
    :[d(.3)r(.8)] *+{C(1,2,1)}="121"
    :"1111"
    )
   ,
   :[r] *+{C(2,2)}="22"  ^{-}
    (
    :"112"  _(.3){-}
    ,
    :[d(.3)r(.8)] *+{C(2,1,1)}="211"
    :"1111"  _{-}
    )
   ,
   :[d(.3)r(.8)] *+{C(3,1)}="31"  _{-}
    (
    :"211"  _(.45){-}
    ,
    :"121"  ^(.7){-}
    )
   )
  ,
  :[u(.3)r(1.2)] *+{C(1,2)}="12"  ^{-}
  :[d(.3)r(.8)] *+{C(1,1,1)}="111"
  ,
  :[d(.3)r(.8)] *+{C(2,1)}="21"  _{-} 
  :"111"  _{-}
  )
 ,
 :[r] *+{C(1,1)}="11" _{-}
 )
)
}$$
Here signs are those coming from the formula \eqref{chdi} for the cochain differential on $E^*(\SW(A_*))$. 
\nl
The following show how quickly the complexity of the cohomology $H^*(E^*_p)$ grows with $p$.
\bex
The cohomology of the horizontal complex of depth $1$ are
$$H^0(E^*_1) = 0\ ,\qquad H^1(E^*_1) = \frac{C(1,1)}{C(2)} \ .$$
The cohomology of the horizontal complex of depth $2$ are
$$H^0(E^*_2) = 0\ ,\qquad H^1(E^*_2) = \frac{C(2,1) \cap C(1,2)}{C(3)}\ ,\qquad H^2(E^*_2) = \frac{C(1,1,1)}{C(2,1) + C(1,2)} \ .$$
The cohomology of the horizontal complex of depth $3$ are
$$H^0(E^*_3) = 0\ ,\qquad H^1(E^*_3) = \frac{C(3,1) \cap C(2,2) \cap C(1,3)}{C(4)}\ ,\qquad H^3(E^*_3) = \frac{C(1,1,1,1)}{C(2,1,1) + C(1,2,1) + C(1,1,2)} \ ,$$
with the degree 2 cohomology being an extension
$$\xymatrix{\frac{C(2,1,1)\ \cap\ C(1,1,2)}{C(2,2)\ +\ C(3,1)\ \cap\ C(1,3)} \ar[r] &\ H^2(E^*_3)\ \ar[r] &\ \frac{C(1,2,1)\ \cap\ \big(C(2,1,1)\ +\ C(1,1,2)\big)}{C(3,1)\ +\ C(1,3)}}\quad .$$
\eex
\bre\lb{ctd}
Note that the highest degree (degree $p$) cohomology of the horizontal complex of an arbitrary depth $p$ is
\beq\lb{hdc}H^p(E^*_p) = \frac{C(1,...,1)}{C(2,1,...,1) + C(1,2,1,...,1)+...+C(1,...,1,2))} \ .\eeq
\ere

In general the degree $m$ component of the horizontal complex of depth $p$ has the form
$$E^m_p = \bigoplus_{|\lambda|=m}C(\lambda)\ ,$$
where $C(\lambda) = C(\lambda_1,...,\lambda_m)$ and the sum is taken over all compositions of $p+1$ into $m$ parts.
\nl
By a {\em composition} $\lambda$ we mean an expansion of $p$ as the sum of a sequence of positive integers $p = \lambda_1+...+\lambda_m$. 
Every composition of $p$ can be thought of as a set-theoretic partition of the set $\{1,...,p\}$ into intervals, i.e. parts of the form $\{\lambda_1+...+\lambda_i,\lambda_1+...+\lambda_i+1,...,\lambda_1+...+\lambda_{i+1}-1\}$. 
Note also that every composition $\lambda$ of $p$ corresponds to a binary vector $x(\lambda)=(x_1,...,x_{p-1})\in\{0,1\}^{p-1}$ of length $p-1$. 
Under this correspondence $x_j$ is $1$ if $j$ and $j+1$ are in different parts of $\lambda$ and $0$ otherwise. 
This correspondence shows that the combinatorial shape of the horizontal complex of depth $p$ is the one of the $(p-1)$-dimensional cube: the centraliser $C(\lambda)$ sits in the vertex with the coordinates $x(\lambda)$, while the edges of the cube correspond to the components of the differential.  Note that the sign with which a differential corresponding to the cube edge $x\to x+e_i$ appears in the chain differential \eqref{chdi} is $(-1)^{\sum_{s=1}^ix_i}$.

Define a {\em cubic diagram} as a collection of subspaces $Q_x$ one 
for each $x\in\{0.1\}^n$ such that for any $x,y\in\{0,1\}^n$ 
$$Q_{xy}\subset Q_x\cap Q_y\ .$$
Here $xy$ is the component-wise product of binary vectors. 
We will also use the labelling by composition $Q(\lambda)$ such that $Q(\lambda) = Q_{x(\lambda)}$.
\nl
A cubic daigram $Q_*$ gives rise to a cochain complex $C^*= C^*(Q_*)$ with
$$C^k = \bigoplus_{l(\lambda)=k-1}Q(\lambda) = \bigoplus_{\sum_ix_i=k}Q_x\ ,$$
where the first sum is taken over all compositions in to $k-1$ parts, while the second sum is over all vectors $x\in\{0.1\}^n$ with the coordinate sum being $k$. 
The differential $d:C^k\to C^{k+1}$ is the direct sum of signed embeddings $Q(\lambda)\to Q(\mu)$, where $\mu$ is the result of subdivision of one of the parts of $\lambda$ into two. Alternatively the differential $d:C^k\to C^{k+1}$ is the direct sum of signed embeddings $Q_x\subset Q_{x+e_i}$, whenever $x_i=0$.
The sign of the embedding $Q(\lambda)\to Q(\mu)$ is $(-1)^j$ if $\mu$ is the division of the $j$-th part of $\lambda$, while the sign of the embedding $Q_x\subset Q_{x+e_i}$ is $(-1)^{\sum_{s=1}^ix_i}$.

\bex
The horizontal complex of depth 3 has the following cubic arrangement
$$
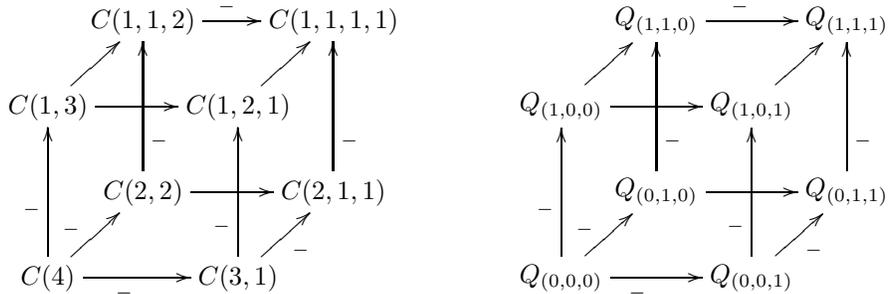
\xygraph{!{0;/r3pc/:;/u2.7pc/::} 
!{(0,2)}*+{C(1,3)}="tl"  !{(2,2)}*+{C(1,2,1)}="tr"
!{(1,1)}*+{C(2,2)}="ml"  !{(3,1)}*+{C(2,1,1)}="mr"
!{(0,0)}*+{C(4)}="bl"  !{(2,0)}*+{C(3,1)}="br"
!{(1,3)}*+{C(1,1,2)}="lt"  !{(3,3)}*+{C(1,1,1,1)}="rt"
"bl":"br" _(.4){-}
"bl":"tl"  ^(.4){-}
"br":"tr"  ^(.3){-}
"tl":"tr"
"bl":"ml"  ^(.4){-}
"ml":"mr"  "br":"mr"  _{-}   "mr":"rt"  _(.3){-}
"ml":"lt" _(.3){-}
 "lt":"rt"   ^(.44){-}
 "tl":"lt"  "tr":"rt"
}\qquad\qquad
\xygraph{!{0;/r3pc/:;/u2.7pc/::} 
!{(0,2)}*+{Q_{(1,0,0)}}="tl"  !{(2,2)}*+{Q_{(1,0,1)}}="tr"
!{(1,1)}*+{Q_{(0,1,0)}}="ml"  !{(3,1)}*+{Q_{(0,1,1)}}="mr"
!{(0,0)}*+{Q_{(0,0,0)}}="bl"  !{(2,0)}*+{Q_{(0,0,1)}}="br"
!{(1,3)}*+{Q_{(1,1,0)}}="lt"  !{(3,3)}*+{Q_{(1,1,1)}}="rt"
"bl":"br" _(.4){-}
"bl":"tl"  ^(.4){-}
"br":"tr"  ^(.3){-}
"tl":"tr"
"bl":"ml"  ^(.4){-}
"ml":"mr"  "br":"mr"  _{-}   "mr":"rt"  _(.3){-}
"ml":"lt" _(.3){-}
 "lt":"rt"   ^(.44){-}
 "tl":"lt"  "tr":"rt"
}
$$
\eex

Now we give a  conditions on the multiplicative sequence $A_*$, which somewhat simplifies the computation of the deformation cohomology $H^n(\SW(A_*))$.
\nl
We say that a multiplicative sequence $A_*$ is {\em generated by its first two members $A_1$ and $A_2$} if the images of the homomorphisms 
$$\mu_{1,...,1,2,1,...,1}:A_1^{\ot m}\ot A_2\ot A_1^{\ot n-m-2}\to A_n$$
jointly generate $A_n$ as an algebra for any $n$.
\nl
Note that in this case the centralisers \eqref{cen} satisfy the property
\beq\lb{jce}C(\lambda)\cap C(\mu) = C(\lambda\cup\mu)\ .\eeq 
Here $\lambda$ and $\mu$ are compositions of $n$ and $\lambda\cup\mu$ stands for their union (the composition corresponding to the union of the equivalence relations corresponding to $\lambda$ and $\mu$). For example $C(2,1) \cap C(1,2) = C(3)$. 

\void{
We also say that the centralisers \eqref{cen} obey the {\em distributivity property} if
\beq\lb{cpr}C(\lambda)\cap\big(C(\mu_1)+...+C(\mu_m)\big) = C(\lambda)\cap C(\mu_1)+...+C(\lambda)\cap C(\mu_m)\ .\eeq

\ble\lb{chc}
Let $A_*$ be a multiplicative sequence generated by its first two members $A_1$ and $A_2$ and such that its centralisers satisfy the distributivity property.
Then the horizontal complex $E^*_p$ of depth $p$ is acyclic away from degree $p$.
\ele
\bpf

\epf

Now we are ready to state our main technical result. 
\bth\lb{dcg}
Let $A_*$ be a multiplicative sequence generated by its first two members $A_1$ and $A_2$ and such that its centralisers satisfy the distributivity property.
Then the deformation cohomology of the Schur-Weyl category $\SW(A_*)$ coincides with the cohomology of the complex
\beq\lb{ccdc}\xymatrix{...\ar[rr] && H^p(E^*_p) \ar[rrr]^(.45){\mu_{1,p} + (-1)^{p+1}\mu_{p,1}} &&& H^{p+1}(E_{p+1}^*) \ar[rr] && ...}\eeq
Here $H^p(E^*_p)$ is the highest degree cohomology of the horizontal complex of depth $p$ given by \eqref{hdc}. 
\eth
\bpf
According to lemma \ref{chc} the only non-zero terms on the first sheet of the spectral sequence \eqref{ss} are 
$$E^{p,0}_1 = H^p(E_p^*)\ .$$
The differential of the first sheet $d^{p,0}_1:E^{p,0}_1\to E^{p+1,0}_1$ coincides with the map $H^p(E_p^*)\to H^{p+1}(E_{p+1}^*)$ induced by 
$\mu_{1,p} + (-1)^{p+1}\mu_{p,1}:A_p\to A_{p+1}$.
\nl
The spectral sequence degenerates at the second sheet giving the result.
\epf
}

\bex\lb{prim}
The first deformation cohomology of a Schur-Weyl category $\SW(A_*)$ is 
$$H^1(\SW(A_*)) = \{a\in Z(A_1)|\ a\ot 1\ot ...\ot 1 + 1\ot a\ot 1\ot ...\ot 1 + ... + 1\ot ...\ot 1\ot a\in Z(A_n)\quad \forall n\}\ .$$
If $A_*$ is a multiplicative sequence generated by its first two members $A_1$ and $A_2$ then
$$H^1(\SW(A_*)) = \{a\in Z(A_1)|\ a\ot 1+1\ot a\in Z(A_2)\ \}\ .$$
\eex

Finally we give an answer for the deformation cohomology under some very restrictive conditions, which will applicable later.
\bth\lb{dcg}
Let $A_*$ be a multiplicative sequence such that its horizontal complexes are acyclic away from the top degree.
Then the deformation cohomology of the Schur-Weyl category $\SW(A_*)$ coincides with the cohomology of the complex
\beq\lb{ccdc}\xymatrix{...\ar[rr] && H^p(E^*_p) \ar[rrr]^(.45){\mu_{1,p} + (-1)^{p+1}\mu_{p,1}} &&& H^{p+1}(E_{p+1}^*) \ar[rr] && ...}\eeq
\eth
\bpf
The zero sheet of the spectral sequence associated to the filtration \eqref{fil} is 
$$E^{p,q}_0 = F^pE^{p+q}/F^{p+1}E^{p+q} = E_p^{p+q}\ .$$
By the assumption its cohomology is only non-zero at 
$$E^{p,0}_1 = H^p(E_p^*)\ .$$
The differential of the first sheet $d^{p,0}_1:E^{p,0}_1\to E^{p+1,0}_1$ coincides with the map $H^p(E_p^*)\to H^{p+1}(E_{p+1}^*)$ induced by 
$\mu_{1,p} + (-1)^{p+1}\mu_{p,1}\to A_{p+1}$.
\nl
The spectral sequence degenerates at the second sheet giving the result.
\epf

\section{Free symmetric categories}

Here we compute the deformation cohomology of free symmetric categories generated by one object.

\subsection{Cubic diagrams of invariants}\lb{cdi}

We denote by $S_n$ the group of permutations on $n$-symbols, i.e. the $n$-th symmetric group.  
Here we look at cubic diagrams of invariants of symmetric group action.

Let $M$ be a linear representation of a symmetric group $S_n$. 
For a composition $n_1+n_2+...+n_r = n$ define $Q((n_1,n_2,...,n_r)$ to be the subspace of invariants $M^{S_{n_1}\times...\times S_{n_r}}$ with respect to the subgroup $S_{n_1}\times...\times S_{n_r}\subset S_n$. 
It is straightforward that the collection of subspaces $Q((n_1,n_2,...,n_r)$ forms a cubic diagram $Q_*(M)$. 
\bex\lb{rer}
Let $k(S_n)$ be the {\em regular $S_n$-representation}, i.e.  the space of function $f:S_n\to k$ with the $S_n$-action given by 
$$(\sigma.f)(x) = f(x\sigma),\qquad \sigma, x\in S_n\ .$$
The components of the cubic diagram $Q_*(k(S_n))$ have the form 
\beq\lb{ccc}Q((n_1,n_2,...,n_r) = k(S_n)^{S_{n_1}\times...\times S_{n_r}} = k(S_n/S_{n_1}\times...\times S_{n_r})\ .\eeq
Note that the right $S_n$-action on $k(S_n)$ 
$$(f.\tau)(x) = f(\tau x),\qquad \tau, x\in S_n$$
gives rise to an $S_n$-action on the cubic diagram $Q_*(k(S_n))$.
\nl
The cubic diagram $Q_*(k(S_n))$ is special in the sense that for any linear representation $M$ of the symmetric group $S_n$
$$Q_*(M) = (Q_*(k(S_n))\ot M)^{S_n}\ ,$$
where the $S_n$-invariants are taken with respect to the diagonal action.
\eex

We will (following the idea of \cite[proposition 2.2]{dr}) compute the cohomology of such cubic diagrams by relating them to certain standard (co)simplicial complexes.
We start by recalling basic facts about simplicial sets.

Let $\Del$ be the category of finite linearly ordered sets and order preserving maps. Denote by $[n]$  the  linearly ordered set $\{0,1,\ldots,n\}$. Any order preserving map $[n]\to[m]$ is a composite of $d_i:[k-1]\rightarrow [k]$ (where $d_i$ is the unique injective map that does not take the value $i\in [k]$) and $s_j:[k]\rightarrow[k-1]$ (where  $s_j$ is the unique surjective map that takes the value $j$ twice). 
\nl
A \textit{simplicial set} is a functor $X:\Del^{op}\rightarrow\Set$. 
We refer to elements of $X(k) = X([k])$ as {\em $k$-simplices} of $X$.
The maps $X(d_i):X(k)\to X(k-1)$ and $X(s_j):X(k-1)\to X(k)$ are called the {\em face} and {\em degeneration maps} correspondingly. 
A simplex is {\em degenerate} if it is in the image of a degeneration map.

\bex
The {\em simplicial interval} is the simplicial set $I=\Del(-,[1]):\Del^{op}\rightarrow\Set$. 
Its $k$-simplices $I(k) = \Del([k],[1])$ are order preserving maps $[k]\to[1]$, which will be represented by non-decreasing binary vectors of length $k+1$. In this presentation the face map $I(d_i)$ removes the $i$-th coordinate, while the degeneration map $I(s_j)$ duplicates the $j$-th coordinate. In particular,  a simplex is degenerate if the corresponding binary vector has repeating coordinates. 
Here is the picture of non-degenerate simplices of the simplicial interval
$$\xymatrix{(0) \ar@{-}[rr]^{(01)} && (1)}$$
\eex

The {\em product} $X \times Y$ of two simplicial sets $X,Y: \Del^{op}\rightarrow \Set$ is the simplicial set 
$$(X \times Y)(k) = X(k) \times Y(k)\ ,\qquad (X \times Y)(f) = X(f) \times Y(f)\qquad f\in \Del([n],[m])\ .$$

\bex\lb{scu}
The {\em simplicial $n$-cube} $I^n$ is the $n$-fold product  $I^{\times n}$. 
We record its simplices as vertical arrays of simplices of $I$, i.e. $k$-simplices $I^n(k)$ are represented by binary matrices with $n$ rows and $k+1$ columns, with non-decreasing rows.
The face map $I^n(d_i)$ removes the $i$-th column, while the degeneration map $I^n(s_j)$ duplicates the $j$-th column. A simplex of $I^n$  is degenerate if the corresponding binary matrix has repeating columns. 
Here is the picture of non-degenerate simplices of the simplicial square $I^2$
$$\xygraph{!{0;/r5pc/:;/u4.7pc/::} 
!{(0,2)}*+{\text{\scalebox{0.6}{$\left(\begin{array}{c}0 \\ 1\end{array}\right)$}}}="tl"  
!{(2,2)}*+{\text{\scalebox{0.6}{$\left(\begin{array}{c}1 \\ 1\end{array}\right)$}}}="tr"
!{(0,0)}*+{\text{\scalebox{0.6}{$\left(\begin{array}{c}0 \\ 0\end{array}\right)$}}}="bl"  
!{(2,0)}*+{\text{\scalebox{0.6}{$\left(\begin{array}{c}1 \\ 0\end{array}\right)$}}}="br"
!{(1,-.2)}*+{\text{\scalebox{0.6}{$\left(\begin{array}{cc}0 & 1\\ 0& 0\end{array}\right)$}}}
!{(1,2.2)}*+{\text{\scalebox{0.6}{$\left(\begin{array}{cc}0 & 1\\ 1& 1\end{array}\right)$}}}
!{(-.3,1)}*+{\text{\scalebox{0.6}{$\left(\begin{array}{cc}0 & 0\\ 0& 1\end{array}\right)$}}}
!{(2.3,1)}*+{\text{\scalebox{0.6}{$\left(\begin{array}{cc}1 & 1\\ 0& 1\end{array}\right)$}}}
!{(1.65,1.3)}*+{\text{\scalebox{0.6}{$\left(\begin{array}{cc}0 & 1\\ 0& 1\end{array}\right)$}}}
!{(1.4,.6)}*+{\text{\scalebox{0.6}{$\left(\begin{array}{ccc}0 & 1& 1\\ 0& 0& 1\end{array}\right)$}}}
!{(.6,1.4)}*+{\text{\scalebox{0.6}{$\left(\begin{array}{ccc}0 & 0& 1\\ 0& 1& 1\end{array}\right)$}}}
"bl"-"br" 
"bl"-"tl"  
"br"-"tr"  
"tl"-"tr"
"bl"-"tr" 
}$$
\eex

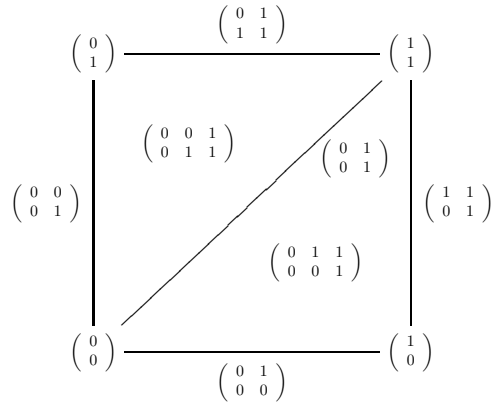
\bex
The {\em boundary of the simplicial $n$-cube} $\partial I^n$ is the simplicial subset of $I^n$ consisting of simplices, which do not have the diagonal 
$$\left(\begin{array}{cc}0 & 1\\ \vdots & \vdots\\ 0& 1\end{array}\right)$$
as an iterated face. In other words, $k$-simplices $\partial I^n(k)$ are represented by binary matrices with either the first column being not all zeroes or the last column being not all ones.
Here is the picture of non-degenerate simplices of the boundary of the simplicial square $\partial I^2$
$$\xygraph{!{0;/r4pc/:;/u3.7pc/::} 
!{(0,2)}*+{\text{\scalebox{0.6}{$\left(\begin{array}{c}0 \\ 1\end{array}\right)$}}}="tl"  
!{(2,2)}*+{\text{\scalebox{0.6}{$\left(\begin{array}{c}1 \\ 1\end{array}\right)$}}}="tr"
!{(0,0)}*+{\text{\scalebox{0.6}{$\left(\begin{array}{c}0 \\ 0\end{array}\right)$}}}="bl"  
!{(2,0)}*+{\text{\scalebox{0.6}{$\left(\begin{array}{c}1 \\ 0\end{array}\right)$}}}="br"
!{(1,-.25)}*+{\text{\scalebox{0.6}{$\left(\begin{array}{cc}0 & 1\\ 0& 0\end{array}\right)$}}}
!{(1,2.25)}*+{\text{\scalebox{0.6}{$\left(\begin{array}{cc}0 & 1\\ 1& 1\end{array}\right)$}}}
!{(-.35,1)}*+{\text{\scalebox{0.6}{$\left(\begin{array}{cc}0 & 0\\ 0& 1\end{array}\right)$}}}
!{(2.35,1)}*+{\text{\scalebox{0.6}{$\left(\begin{array}{cc}1 & 1\\ 0& 1\end{array}\right)$}}}
"bl"-"br" 
"bl"-"tl"  
"br"-"tr"  
"tl"-"tr"
}$$
\eex

For a simplicial set $X$ denote by $C^*(X)$ its {\em normalised cochain complex} with coefficients $k$, i.e. the collection of 
$$C^n(X) = \{f\in \Set(X(n),k)|\ f\circ X(s_j) = 0\quad \forall j\}$$ 
vector spaces of functions on non-degenerate $n$-simplices of $X$ with the differential 
$$\partial:C^n(X)\to C^{n+1}(X),\qquad \partial = \sum_{i=0}^n(-1)^i\partial_i\ ,$$
where $\partial_i(f) = f\circ X(d_i)$. 
\nl
Note that a map of simplicial sets $X\to Y$ (i.e. a natural transformation of functors $X,Y: \Del^{op}\rightarrow \Set$) gives rise to 
a homomorphism of cochain complexes $C^*(Y)\to C^*(X)$.
\nl
For a simplicial subset $X\subset Y$ denote by $C^*(Y,X) = Ker(C^*(Y)\to C^*(X))$ its {\em relative normalised cochain complex}.

Note that the natural $S_n$-action on the simplicial $n$-cube $I^n$ gives rise to an $S_n$-action on the relative normalised cochain complex $C^*(I^n,\partial I^n)$.
\ble\lb{rcc}
The relative normalised cochain complex $C^*(I^n,\partial I^n)$ is isomorphic in an $S_n$-equivariant way to the cochain complex $C^*(Q_*(k(S_n)))[1]$ of the cubic diagram of the regular representation of $S_n$ shifted by one in degree.
\ele
\bpf
The cochain complex $C^*(I^n,\partial I^n)$ in degree $m$ is the space of functions on non-degenerate $m$-simplices of $I^n$, which contain the diagonal. In the notations of example \ref{scu} these simplices correspond to $n$-by-$(m+1)$ binary matrices with non-decreasing rows, with the first column being all zeroes and with the last column being all ones.  
It is straightforward to see that by a permutation of rows (which corresponds to the $S_n$-action on $C^*(I^n,\partial I^n)$) any such matrix can be brought to a form
\beq\lb{cfm}\left(\begin{array}{ccccc} 
0 & 0 & \quad \hdots \quad & 0 & 1\\
\vdots&\vdots&& \vdots & \vdots\\
\vdots &\vdots & & 0 & \vdots\\
\vdots &\vdots && 1 & \vdots \\
\vdots&\vdots & \text{\reflectbox{$\ddots$}} & \vdots & \vdots\\
\vdots &0 && \vdots & \vdots \\
\vdots &1 && \vdots & \vdots \\
\vdots &\vdots && \vdots & \vdots \\
0 &1 &\hdots & 1 & 1 \\
\end{array}\right)\eeq
where the number of zeroes in each row is not more than in the previous one.
Let $n_{m-k}$ be the increment of the number of ones in the $(k+2)$-th column compared to $(k+1)$-th column. Clearly $n_1+n_2+...+n_r = n$ is a composition and $S_{n_1}\times...\times S_{n_r}\subset S_n$ is the $S_n$-stabiliser of the matrix \eqref{cfm}. 
Thus the subspace $C^m(I^n,\partial I^n)$ of functions supported on the $S_n$-orbit of \eqref{cfm} coincides with the component \eqref{ccc} of the cubic diagram $Q_*(k(S_n))$. 
\nl
The differential of  $C^*(I^n,\partial I^n)$ is the alternating sum of maps of functions induced by column erasing maps. It is easy to see that they correspond to the embeddings of the cubic diagram $Q_*(k(S_n))$.
\epf

Here is our main technical tool for computing the deformation cohomology of free symmetric categories.
\bpr\lb{chc}
Let $M$ be a linear representation of a symmetric group $S_n$. 
The cohomology of the cochain complex $C^*(Q_*(M))$ of the cubic diagram $Q_*(M)$ of invariants is zero apart from degree $n-1$, where it has the form
$$H^{n-1}(Q_*(M)) = \frac{M}{\sum_{i=1}^{n-1}(1+t_i)(M)}\ .$$
\epr
\bpf
By example \ref{rer} the cochain complex $C^*(Q_*(M))$ coincides with $(Q_*(k(S_n))\ot M)^{S_n}$.
By lemma \ref{rcc} the cochain complex $C^*(Q_*(k(S_n)))$ is isomorphic to the degree shift of the relative normalised cochain complex $C^*(I^n,\partial I^n)[-1]$.
Since the geometric realisation of the simplicial cube $I^n$ is homeomorphic to the $n$-disk the cohomology $H^*(I^n)$ is concentrated in degree zero. Similarly,  the geometric realisation of the boundary $\partial I^n$ of the simplicial cube is homeomorphic to the $(n-1)$-sphere and  its the cohomology $H^*(I^n)$ is non-zero only in degree zero and $n-1$.
The long exact sequence of the relative cohomology shows that the cohomology $H^*(I^n,\partial I^n)$ is non-zero only in degree $n$.
Finally remark \ref{ctd} gives the answer for the top cohomology $H^{n-1}(Q_*(M))$. 
\epf

\subsection{The free symmetric category on one object}

Here we assume that the characteristic of the ground field $k$ is odd.
Recall from  \cite{dm} that the Schur-Weyl category $\SW(A_*)$ of the multiplicative sequence $A_n = k[S_n]$ is the free symmetric category $\cS$ on one object $X$.
In this case the centralisers $C(\lambda)$ can be written as the subalgebras of invariants $k[S_n]^{S_{\lambda_1}\times...\times S_{\lambda_m}}$ with respect to the conjugation action. 
In other words the horizontal complexes correspond to cubic diagrams of $S_n$-invariant of the adjoint action on $k[S_n]$.  
According to proposition \ref{chc} the horizontal complex of depth $n$ is acyclic away from the top degree $n$ and its $n$-th cohomology is 
\beq\lb{tpfs}H^n(E^*_n) = k[S_n]/\left(C_{k[S_n]}(t_1) + C_{k[S_n]}(t_2) + ... + C_{k[S_n]}(t_{n-1})\right) \ .\eeq
Assuming that the ground field $k$ is of odd characteristic each centraliser $C_{k[S_n]}(t_i)$ can be written as $\{x + t_ixt_i |\ x\in k[S_n]\}$.
Thus 
a linear function on $H^n(E^*_n)$ is a function $f:S_n\to k$, which is zero on $x + t_ixt_i$ for all $x\in k[S_n]$ and all $i=1,...,n-1$.
In other words, the dual space of $H^n(E^*_n)$ is the vector space of functions $f:S_n\to k$ such that 
\beq\lb{sc}f(\sigma\pi\sigma^{-1}) = sign(\sigma)f(\pi)\qquad \forall\ \sigma,\pi\in S_n\eeq
Denote it $V_n$. 
\bpr
Let $k$ be a field of odd characteristic.
Then the deformation cohomology of the free symmetric category $\cS$ is 
$$H^n(\cS) = k[S_n]/\left(C_{k[S_n]}(t_1) + C_{k[S_n]}(t_2) + ... + C_{k[S_n]}(t_{n-1})\right)  \ .$$
\epr
\bpf
By theorem \ref{dcg} the deformation cohomology $H^*(\cS)$ coincides with the cohomology of the complex \eqref{ccdc}.
Its dual complex (according to proposition \ref{scs}) is 
$$\xymatrix{... && V_p \ar[ll] &&& V_{p+1} \ar[lll]_{\mu^*_{1,p} + (-1)^{p+1}\mu^*_{p,1}} && ... \ar[ll]}$$
Note that for any $\sigma\in S_p$ one has $\mu_{1,p}(\sigma) = \tau\mu_{p,1}(\sigma)\tau^{-1}$, where $\tau = (1,2,...,p+1)\in S_{p+1}$ is a $(p+1)$-cycle. Since $sign(\tau) = (-1)^p$ the property \eqref{sc} implies that 
$$\mu^*_{1,p}(f) = f\circ\mu_{1,p}  = (-1)^p f\circ\mu_{p,1} = (-1)^p\mu^*_{p,1}(f)$$ for any $f\in V_{p+1}$.
Thus $\mu^*_{1,p} + (-1)^{p+1}\mu^*_{p,1} = 0$ and $H^p(\cS) = H^n(E^*_n)$.
\epf

\ble\lb{scs}
Let $s_n$ be the dimension of $V_n$.
Then 
$$\sum_{n\geq 0} s_nt^n = \prod_{m\geq 1}(1+t^{2m-1})\ .$$
\ele
\bpf
Note that a function $f\in V_n$ must  be supported on a union of conjugacy classes of $S_n$.
Moreover the centraliser of any conjugacy class in the support should consist of just even permutations. 
It is straightforward to see that such classes correspond to partitions of $n$ into distinct odd parts.
The generating function for the numbers of such partitions is clearly $\prod_{m\geq 1}(1+t^{2m-1})$. 
\epf

The first few terms of the generating series are
$$\sum_{n\geq 0} s_nt^n = 1 + t + t^3 + t^4 + t^5 + t^6 + t^7 + 2t^8 + 2t^9 + 2t^{10} + 2t^{11} + 3t^{12} + ... $$

\bex\lb{exa}
For an odd $n$ define the function $f_n\in V_n$ by $f_n(t_1t_2...t_{n-1}) = 1$.
For example, when $n=3$ the function $f_3(t_1t_2) = -f_3(t_2t_1) = 1$ spans $V_3$. 
\nl
For $f\in V_n$ define $c_f = \sum_{\sigma\in S_n}f(\sigma)\sigma\in k[S_n]$. 
Denote $e_{2i-1} = c_{f_{2i-1}}$. For example $e_1 = 1, e_3 = t_1t_2 - t_2t_1$. 
We will also denote by $e_{2i-1}$ the corresponding class in $H^{2i-1}(\cS)$.
\eex

\bth\lb{cfs}
Let $k$ be a field of odd characteristic.
Then the deformation cohomology of the free symmetric category $\cS$ is the exterior algebra 
$$H^*(\cS) = \Lambda(e_1,e_3,e_5,...),\qquad deg(e_{2i-1}) = 2i-1 $$
generated by elements defined in example \ref{exa}.
\eth
\bpf
The cup-product \eqref{cup} defines the homomorphism $\Lambda(e_1,e_3,e_5,...)\to H^*(\cS)$. 
It follows from the proof of lemma \ref{scs} that this homomorphism is surjective.
Finally dimension comparison shows that this is an isomorphism.
\epf

\subsection{Free symmetric categories}

Let $A$ be a commutative algebra.
Denote by $A^{\ot n}*S_n$ the skew group algebra with respect to the permutation $S_n$-action on $A^{\ot n}$.
The Schur-Weyl category $\C=\SW(A_*)$ of the multiplicative sequence 
$A_n = A^{\ot n}*S_n$ is the free symmetric category $\cS(A)$ on one object $X$ with the endomorphism algebra $End_\C(X)=A$ (see \cite{dm} for details).

Denote by $S^{m}(A) = (A^{\ot m})^{S_{m}}$ the $m$-th symmetric power of a vector space $A$. 
\ble \lb{cfs}
Let $A$ be a commutative algebra without zero divisors, which is more than one dimensional. 
Then
$$C(n_1,...,n_r) = S^{n_1}(A)\ot...\ot S^{n_r}(A)\ .$$
\ele
\bpf
We first consider the largest  centralizer $C(1,1,...1) = C_{A^{\ot n}*S_n}(A^{\ot n})$.
An element  $a=\sum a_\sigma  \sigma\in A^{\ot n}*S_n$ belongs to $C(1,1,...1)$ if for any $b\in A^{\ot n}$ we have 
$$0 = \left[ \sum a_\sigma * \sigma, b \right]  = \sum [ a_\sigma * \sigma , b]  = \sum a_\sigma [ \sigma, b ] = \sum a_\sigma (\sigma(b) - b )  \sigma\ .$$
Thus $a_\sigma (\sigma(b) - b) = 0$. 
Note that if $\sigma$ is non-trivial we can find such $b\in A^{\ot n}$ that $\sigma(b) - b$.  Indeed, since $dim(A) \geq 2$, we can choose two linearly independent elements of $c,d\in A$ $c$ and set $b = c \otimes c \otimes c \otimes d \otimes c ... c \otimes c$, where the position of $d$ is where $\sigma$ has acts nontrivially. 
Hence $a_\sigma=0$ for any non-trivial $\sigma$ and so we get that $C(1,1,...1) =  A^{\otimes n}$. 
\nl
Now note that $C(n_1,n_2,...,n_r) = C(1,1,...1)^{S_{n_1} \times ... \times S_{n_r}}$, so
$$C(n_1,n_2,...,n_r) = (A^{\ot n})^{S_{n_1} \times ... \times S_{n_r}} = S^{n_1}(A)\ot...\ot S^{n_r}(A)\ .$$
\epf

Denote by $\Lambda^n(A)$ the exterior power of $A$. 
Thus we recover the following (originally \cite[proposition 2.2]{dr}).
\ble\lb{hcsa}
Let $A$ be a commutative algebra without zero divisors, which is more than one dimensional. 
Then the horizontal complex of $\SW(A_*)$ of depth $n$ is acyclic away from the top degree, where its cohomology is $\Lambda^n(A)$. 
\ele
\bpf
This is a direct consequence of proposition \ref{chc}, since 
$$ \frac{M}{\sum_{i=1}^{n-1}(1+t_i)(M)} = \Lambda^n(A)\ .$$
\epf

\bex
For $a\in A$ define the endomorphism $\psi(a)\in End(Id_{\cS(A)})$ by
$$\psi(a)_{[n]} = a\ot 1\ot...\ot 1 + 1\ot a\ot 1\ot...\ot 1 + ... + 1\ot...\ot 1\ot a\ \in A^{\ot n}\subset End_{\cS(A)}([n])\ .$$
It follows from example \ref{prim} that the assignment $a\mapsto \psi(a)$ gives an isomorphism
$H^1(\cS(A)) = A$. 
\eex

\bth\lb{cfsa}
Let $A$ be a commutative algebra without zero divisors, which is more than one dimensional. Then
the deformation cohomology of the free symmetric category $\cS(A)$ is the exterior algebra of the first cohomology $H^1(\cS(A)) = A$
$$H^*(\cS(A))\ \simeq\ \Lambda^*(A)\ .$$
\eth
\bpf
The cup-product \eqref{cup} defines the homomorphism $\Lambda^*(A) = \Lambda^*(H^1(\cS(A)))\to H^*(\cS(A))$,
which is an isomorphism by lemma \ref{hcsa} and theorem \ref{dcg}. 
\epf

\bre\lb{pfs}
The theorem \ref{cfs} say that the deformation cohomology of the free symmetric category $\cS(A)$ is generated by the space $Prim(\S(A))$ of primitive endomorphisms of the identity functor.
\nl
Explicitly an endomorphism $\alpha(a)$ corresponding to $a\in A$ has the following specialisations
$$\alpha(a)_{[n]} =  a\ot 1\ot...\ot 1 + 1\ot a\ot ... 1 +  ... + 1\ot ... 1\ot a\in A^{\ot n}*S_n\ .$$
\ere

\subsection{Degenerate affine Hecke category}

The {\it degenerate affine Hecke algebra} $\Lambda_n$ is the
unital associative algebra generated by elements
$t_1,\dots,t_{n-1}$ and $y_1,\dots,y_n$
subject to the relations
\begin{alignat}{2}
t_i^2=1,\qquad t_it_{i+1}t_i&=t_{i+1}t_it_{i+1},
\qquad &&t_it_j=t_jt_i\quad\text{for}\quad |i-j|>1,
\\
y_it_i-t_iy_{i+1}&=1,\qquad &&y_iy_j=y_jy_i.
\end{alignat}
The assignments
\beq
\begin{aligned}
t_i&\otimes 1\mapsto t_i,\qquad 1\otimes t_j\mapsto t_{j+m},\\
y_i&\otimes 1\mapsto y_i,\qquad 1\otimes y_j\mapsto y_{j+m}
\end{aligned}
\eeq
define algebra homomorphisms
\beq
\Lambda_m\otimes\Lambda_n\to\Lambda_{m+n}.
\eeq
It is easy to see that
these homomorphisms satisfy the associativity axiom
thus giving rise to the multiplicative sequence of algebras
$\Lambda_*=\{\Lambda_n\ |\ n\geqslant 0\}$
The Schur-Weyl category $\C(\Lambda_*)$ of the multiplicative sequence
$\Lambda_*$ was called the {\em degenerate affine Hecke category} $\caL$ in \cite{dm}.

Define the {\em length} $l$ of a permutation $\sigma$ as the smallest number of neighbouring transpositions $t_i$ 
required to write $\sigma$ as a product. We extend the length function to elements of $k[S_n]$ by taking the maximum over all non-zero summands:  
$l(\sum a_{\tau}\tau) = \underset{\tau: a_\tau\not=0}{max}\ l(\tau)$. 

For $i=1,...,n$ define $\partial_i:S_n\to \Lambda_n$ by $\partial_i(\sigma) = x_i \sigma - \sigma x_{\sigma^{-1}(i)}$. 
\ble\lb{lda}
The map $\partial_i$ obeys the following  twisted Leibniz formula 
$$ \partial_i ( \sigma \tau)  = \partial_i( \sigma)  \tau + \sigma \partial_{\sigma^{-1}(i)}( \tau) $$ 
Moreover, $\partial_i(\sigma) \in k[S_n]$
and $l(\partial_i (\sigma)) < l(\sigma)$ in other terms, $\partial_i(\sigma)$ is a linear combination of permutations that are less in length than $\sigma$.
\ele  
\bpf 
We first prove the Leibniz formula: 
$$ x_i \sigma \tau =( \sigma x_{\sigma^{-1}(i)} + \partial_i (\sigma))\tau
= \sigma ( \tau x_{\tau^{-1}(\sigma^{-1}(i))}  + \partial_{\sigma^{-1}(i)}(\tau) ) + \partial_i (\sigma) \tau =
$$
$$
 = \sigma \tau x_{    \tau^{-1}({    \sigma^{-1}(i))    } } +  \sigma  \partial_{\sigma^{-1}(i)}(\tau) +  \partial_i (\sigma)\tau\ .$$ 
Hence, 
$$\partial_i ( \sigma \tau ) = x_i \sigma \tau  - \sigma \tau x_{(\sigma\tau)^{-1}(i)}  = \sigma  \partial_{\sigma^{-1}(i)}(\tau) +  \partial_i( \sigma)(\tau)
 $$  
We prove that $\partial_i(\sigma) \in k[S_n]$ by induction on length. 
In the base case, when $l(\sigma)=1$ we have $\partial_i(t_j) = 0$ for $|i-j|>1$, and $\partial_i(t_i)  = 1$, $\partial_{i+1}(t_i)  = -1$, all of which belong to $k[S_n]$. 
Now assume that  $\partial_j(\sigma)\in k[S_n]$ for all $j$ and $\sigma$ where $l(\sigma)<l$. For $l(\sigma')=l$ write $\sigma' = t_j\sigma$ for $\sigma$, where $l(\sigma)=l-1$. By the Leibniz formula $\partial_i(\sigma')= \partial_i(t_i\sigma) = \partial_i(t_j)\sigma + t_j \partial_{t_j^{-1}(i)}(\sigma)$, which clearly belongs to $k[S_n]$. 
\nl
We also prove the decreasing length property by induction on length.
Again in the base case
$\partial_i(t_j) = 0$ for $|i-j|>1$, and $\partial_i(t_i)  = 1 $, $\partial_{i+1}(t_i)  = -1 $, which fits the property $l(\partial_i (\sigma)) < l(\sigma)$ since scalars have zero  length. 
Again taking $l(\sigma')=l$ and writing it as $\sigma' = t_j\sigma$ with $l(\sigma)=l-1$ we get
$\partial_i(\sigma') = \partial_i(t_j \sigma) = \partial_i(t_j)\sigma + t_j \partial_{t_j^{-1}(i)}(\sigma)$.
Since $l(\sigma) < l(t_j \sigma)$ and $\partial_i(t_j)$ is just a scalar as seen above, the first term has length less than $t_j\sigma$. By induction for the second term we have $l(\partial_{t_j^{-1}(i)}(\sigma)) < l(\sigma)$, hence $l(t_j \partial_{t_j^{-1}(i)}(\sigma)) < l(t_j \sigma)$. Thus the second term also satisfies the property and we are done.
\epf 

The following is a generalisation of the Bernstein's lemma, computing the centre of a degenerate affine Hecke algebra.
\ble 
The centralisers $C(m_1,...,m_n) = C_{\Lambda_{m_1+...+m_n}}(\Lambda_{m_1}\ot...\ot \Lambda_{m_n})$ have the form
$$C(m_1,...,m_n) = S^{m_1}(k[x])\ot...\ot S^{m_n}(k[x])\ .$$
\ele
\bpf
The proof is very similar to the one of lemma \ref{cfs}. 
We first show that $C(1,1,...,1) = k[x]^{\ot n}$. 
Note that an element $a\in \Lambda_n$ can be written as a unique combination $a = \sum_{\sigma\in S_n}  a_\sigma \sigma$ with $a_\sigma \in k[x]^{\ot n}$. Note also that $a$ belongs to $C(1,1,..1)$ iff the commutator $[a,x_i]$ is zero.  Writing
$$[a,x_i] = \sum_\sigma a_\sigma x_i \sigma - \sum_\sigma a_\sigma (x_{\sigma(i)} \sigma + \partial_{\sigma(i)}(\sigma) ) = \sum_\sigma a_\sigma (x_i - x_{\sigma(i)}) \sigma  +\sum_\sigma a_\sigma\partial_{\sigma(i)}(\sigma)$$
we get
$$\sum_\sigma  a_\sigma\partial_{\sigma(i)}(\sigma) = \sum_\sigma a_\sigma (x_{\sigma(i)} - x_i) \sigma $$
Take $\sigma_0$ to be a permutation of maximum length such that $a_{\sigma_0} \neq 0$. 
By lemma \ref{lda} we have that $l( \partial_{\sigma(i)} \sigma) < l(\sigma) < l(\sigma_0)$, thus $l(\sum_\sigma a_\sigma (x_i - x_{\sigma(i)}) \sigma) < l(\sigma_0)$, which implies that the coefficient of $\sigma_0$ in $\sum_\sigma a_\sigma (x_i - x_{\sigma(i)}) \sigma$ is zero, i.e. $a_{\sigma_0} (x_i - x_{\sigma_0(i)}) =0$. Since $ a_{\sigma_0} \neq 0$ we have $x_i - x_{\sigma_0(i)} = 0$ for all $i$, which implies $\sigma_0 = e$. Thus $a\in k[x]^{\ot n}$ as desired. 
\nl
To treat the general case of $C(m_1,...,m_n)$ recall a well known formula 
$$ at_i - t_it_i(a)  = \frac{a - t_i(a)}{x_i - x_{i+1}}\ ,$$ 
where $a$ is an element of $k[x]^{\ot n}$. By the definition, an element $a$ of $C(1,1,...,1) = k[x]^{\ot n}$ belongs to $C(m_1,...,m_n)$ iff the commutator $[a,t_i]$ is zero for any $t_i\in S_{m_1} \times ... \times S_{m_n}$. Writing
$$at_i - t_ia = at_i - t_it_i(a) - t_i(a- t_i(a)) = \frac{a - t_i(a)}{x_i - x_{i+1}}- t_i(a- t_i(a))$$
we get $a- t_i(a)=0$. 
Thus 
$$C(m_1,...,m_n) = C(1,1,...,1) ^{S_{m_1} \times ... \times S_{m_n}} = S^{m_1}(k[x])\ot...\ot S^{m_n}(k[x])\ .$$
\epf

\bth\lb{dvda}
The deformation cohomology of the degenerate affine Hecke category $\caL$ is the exterior algebra of $k[x]$
$$H^*(\caL)\ \simeq\ \Lambda^*(k[x])\ .$$
\eth
\bpf
The proof is identical to the proof of theorem \ref{cfs} for $A = k[x]$. 
\epf

\section{Representations of general linear Lie algebras}

Denote by $\Rep(\g)$ the symmetric tensor category of representations of a Lie algebra $\g$. 

\subsection{Deformation cohomology of $\Rep(\g)$}

Here we recall (from \cite{da}) some basic facts about deformation cohomology of categories of modules over Hopf algebras in general and of representations of  Lie algebras in particular.

Let $H$ be a Hopf algebra. 
Denote by $H\da\M od$ the tensor category of its (left) modules.
Denote by $F:H\da\M od\to\Vect$ the forgetful functor. 

The collection of tensor powers $H^{\ot n}$ has a structure of cosimplicial complex with
coface maps ${\partial}_{i}:H^{n}\longrightarrow H^{n+1}$
$${\partial}_{i}(h_{1}\otimes ...\otimes h_{n}) =
\left\{
\begin{array}{cl}
1\otimes h_{1}\otimes ...\otimes h_{n},& i=0\\
h_{1}\otimes ...\otimes\Delta (h_{i})\otimes ...\otimes h_{n},& 1\leq i\leq n\\
h_{1}\otimes ...\otimes h_{n}\otimes 1,& i=n+1
\end{array}
\right.
$$
and codegeneracy maps 
$${\sigma}_{i}(h_{1}\otimes ...\otimes h_{n+1}) = h_{1}\otimes ...\otimes\varepsilon (h_{i})\otimes ...\otimes h_{n+1}))$$

The associated (unnormalised) cochain complex $C^*_{coHoch}(H)$ is called the {\em co-Hochschild complex}
and its cohomology - the {\em co-Hochschild cohomology} of Hopf algebra $H$. 

\bpr\lb{dch}
The deformation complex of the forgetful functor $F:H\da\M od\to\Vect$ is isomorphic to the co-Hochschild complex of $H$.
\nl
The deformation complex of the tensor category $H\da\M od$ is isomorphic to the subcomplex of adjoint $H$-invariants of the co-Hochschild complex of $H$.
\epr

\bre
The isomorphism in the above proposition is realized by two mutually inverse maps:
$$End(F^{\otimes n})\longrightarrow H^{\otimes n}$$
which sends the endomorphism to its specialization on the objects $H,...,H$,
and
$$H^{\otimes n}\longrightarrow End(F^{\otimes n})\ ,$$
which associates to the element $x\in H^{\otimes n}$ the endomorphism of
multiplying by $x$.
\ere

Note that the deformation complexes of the forgetful functor and of the category of modules over a Hopf algebra depend only on the coalgebra structure. 

For a Lie algebra $\g$ the category of representations $\Rep(\g)$ coincides with the category of modules $U(\g)\da\M od$ over the universal enveloping algebra $U(\g)$. 
The following was also established in \cite{da}. We add a sketch of the proof here.
\bth\lb{dcrc}
The deformation cohomology of the forgetful functor $F:\Rep(\g)\to\Vect$ is the exterior algebra of $\g$
$$H^*(F)\ \simeq\ \Lambda^*(\g)\ .$$
\nl
The deformation cohomology of the tensor category of representations of the Lie algebra $\g$ is the subalgebra of invariants
$$H^*(\Rep(\g))\ \simeq\ \Lambda^*(\g)^\g$$
of the exterior algebra of $\g$. 
\eth
\bpf
By Poincare-Birkhoff-Witt theorem $U(\g)$ is isomorphic to the symmetric algebra $S^*(\g)$ as a coalgebra. 
The grading of $S^*(\g)$ induces a grading on the co-Hochschild complex $C^*_{coHoch}(S^*(\g))$. 
The graded component of degree $n$ of $C^*_{coHoch}(S^*(\g))$ is the cubic diagram $Q_*(\g^{\ot n})$ of invariants the natural $S_n$-module $\g^{\ot n}$. The first part of the theorem now follows from proposition \ref{chc}.
\nl
By proposition \ref{dch} the deformation complex $C^*(\Rep(\g))$ is isomorphic to the subcomplex $C^*_{coHoch}(U(\g))^\g$ of adjoint $\g$-invariants.
Since the Poincare-Birkhoff-Witt isomorphism $S^*(\g)\to U(\g)$ is equivariant with respect to the adjoint $\g$-action
the latter is isomorphic to the subcomplex $C^*_{coHoch}(S^*(\g))^\g$. 
By the above the cubic diagram $Q_*(\g^{\ot n})$ is quasi-isomorphic to its cohomology $\Lambda^n(\g)[n]$, with the quasi-isomorphism given by the projection $\g^{\ot n}\to\Lambda^n(\g)$. Since the projection is $\g$-equivariant the subcomplex $Q_*(\g^{\ot n})^\g$ is also quasi-isomorphic to its cohomology $\Lambda^n(\g)^\g[n]$. 
\epf

\subsection{Deformation cohomology of $\Rep(\gl(V))$}

Let $V$ be a vector space over $k$ of dimensional $dim(V)=d$. 
Here $V^*$ stands for the dual vector space.
Denote by $\gl(V)=V\ot V^*$ the general linear Lie algebra of $V$. 

We will use Penrose's calculus \cite{pe} to represent tensors graphically. 
Our diagrams are to be read top down.
We will use orientation to distinguish between $V$ and $V^*$. 
For example, the identity morphisms on $V$ and $V^*$ are 
$$
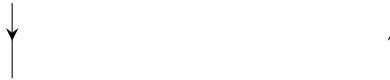
\begin{tikzpicture}
\draw [middleprightarrow] (0,1) to (0,0);
\end{tikzpicture}
\hspace{12pc}
\begin{tikzpicture}
\draw [middlepleftarrow] (0,1) to (0,0);
\end{tikzpicture}
$$
correspondingly.
The canonical pairing $V^*\ot V\to k$ and the canonical element $k\to V\ot V^*$ have the form
$$
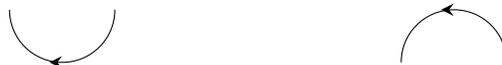
\begin{tikzpicture}
\draw [middlepleftarrow] (0,0.7) to [out=-90,in=180] (.7,0) to [out=0,in=-90] (1.4,0.7);
\end{tikzpicture}
\hspace{9pc}
\begin{tikzpicture}
\draw [middlepleftarrow] (0,0) to [out=90,in=180] (.7,.7) to [out=0,in=90] (1.4,0);
\end{tikzpicture}
$$
The following is the graphical presentation of the adjunction between the canonical pairing and the canonical element  
$$
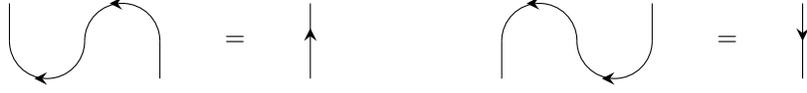
\begin{tikzpicture}
\draw  (0,0) to (0,.5);
\draw [middlepleftarrow] (0,0) to [out=-90,in=180] (.5,-.5) to [out=0,in=-90] (1,0);
\draw [middlepleftarrow] (1,0) to [out=90,in=180] (1.5,.5) to [out=0,in=90] (2,0);
\draw  (2,0) to (2,-.5);
\node at (3,0) {=};
\draw [middlepleftarrow] (4,.5) to (4,-.5);
\end{tikzpicture}
\hspace{6pc}
\begin{tikzpicture}
\draw  (0,0) to (0,-.5);
\draw [middlepleftarrow] (0,0) to [out=90,in=180] (.5,.5) to [out=0,in=90] (1,0);
\draw [middlepleftarrow] (1,0) to [out=-90,in=180] (1.5,-.5) to [out=0,in=-90] (2,0);
\draw  (2,0) to (2,.5);
\node at (3,0) {=};
\draw [middleprightarrow] (4,.5) to (4,-.5);
\end{tikzpicture}
$$
The transposition of tensor factors is depicted as
$$
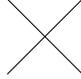
\begin{tikzpicture}
\draw  (0,.5) to (1,-.5);
\draw  (1,.5) to (0,-.5);
\end{tikzpicture}
$$
For example, the picture
$$
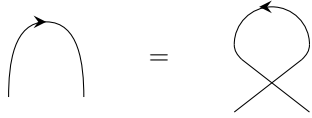
\begin{tikzpicture}
\draw [middleprightarrow] (0,-.5) to [out=90,in=180] (.5,.5) to [out=0,in=90] (1,-.5);
\node at (2,0) {=};
\draw [middlepleftarrow] (3,.2) to [out=90,in=180] (3.5,.7) to [out=0,in=90] (4,.2);
\draw (3,.2) to [out=-90,in=135] (3.1,0);
\draw (3.1,0) to (4,-.7);
\draw (4,.2) to [out=-90,in=45] (3.9,0);
\draw (3.9,0) to (3,-.7);
\end{tikzpicture}
$$
represents the canonical element $k\to V^*\ot V$. 
Another example is the identity
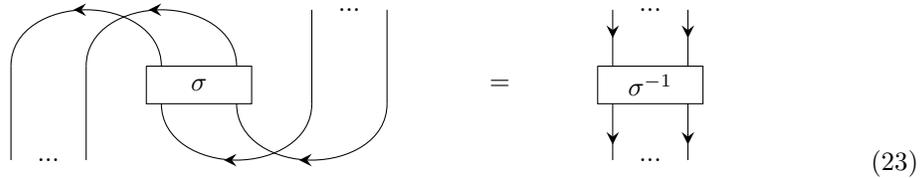
\beq\lb{inp}
\begin{tikzpicture}
\draw (-2,-1) to (-2,.25);
\draw [middlepleftarrow] (-2,.25) to [out=90,in=180] (-1,1) to [out=0,in=90] (0,.25);
\node at (-1.5,-1) {...};
\draw (-1,-1) to (-1,.25);
\draw [middlepleftarrow] (-1,.25) to [out=90,in=180] (0,1) to [out=0,in=90] (1,.25);
\draw  (-.2,.25) to (1.2,.25) to (1.2,-.25) to (-.2,-.25) to (-.2,.25);
\node at (.5,0) {$\sigma$};
\draw [middlepleftarrow] (0,-.25) to [out=-90,in=180] (1,-1) to [out=0,in=-90] (2,-.25);
\draw (2,1) to (2,-.25);
\node at (2.5,1) {...};
\draw [middlepleftarrow] (1,-.25) to [out=-90,in=180] (2,-1) to [out=0,in=-90] (3,-.25);
\draw (3,1) to (3,-.25);
\node at (4.5,0) {=};
\draw [middleprightarrow] (6,1) to (6,.25);
\node at (6.5,1) {...};
\draw [middleprightarrow] (7,1) to (7,.25);
\draw  (5.8,.25) to (7.2,.25) to (7.2,-.25) to (5.8,-.25) to (5.8,.25);
\node at (6.5,0) {$\sigma^{-1}$};
\draw [middlepleftarrow] (6,-1) to (6,-.25);
\node at (6.5,-1) {...};
\draw [middlepleftarrow] (7,-1) to (7,-.25);
\end{tikzpicture}
\eeq
where $\sigma\in S_m$ is a permutation. 

As above $alt_m = (m!)^{-1}\sum_{\sigma\in S_m} sign(\sigma)\sigma$ stands for the projector on the $m$-th exterior power.
\nl
For a permutation $\sigma\in S_m$ denote by $\sigma^{(2)}:(V\ot V^*)^{\ot m}\to (V\ot V^*)^{\ot m}$ the automorphism permuting the factors $V\ot V^*$. Graphically
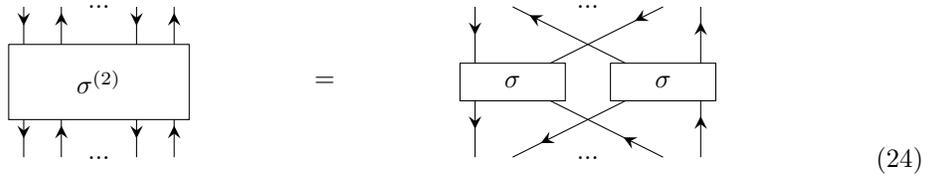
\beq\lb{dup}
\begin{tikzpicture}
\draw [middleprightarrow] (0,1) to (0,0.5);
\draw [middlepleftarrow] (0.5,1) to (0.5,0.5);
\node at (1,1) {...};
\draw [middleprightarrow] (1.5,1) to (1.5,0.5);
\draw [middlepleftarrow] (2,1) to (2,0.5);
\draw  (-.2,.5) to (2.2,.5) to (2.2,-.5) to (-.2,-.5) to (-.2,.5);
\node at (1,0) {$\sigma^{(2)}$};
\draw [middleprightarrow] (0,-.5) to (0,-1);
\draw [middlepleftarrow] (0.5,-.5) to (0.5,-1);
\node at (1,-1) {...};
\draw [middleprightarrow] (1.5,-.5) to (1.5,-1);
\draw [middlepleftarrow] (2,-.5) to (2,-1);
\node at (4,0) {=};
\draw [middleprightarrow] (6,1) to (6,0.25);
\draw [leftpleftarrow] (6.5,1) to (8,0.25);
\node at (7.5,1) {...};
\draw [leftprightarrow] (8.5,1) to (7,0.25);
\draw [middlepleftarrow] (9,1) to (9,0.25);
\draw  (5.8,.25) to (7.2,.25) to (7.2,-.25) to (5.8,-.25) to (5.8,.25);
\node at (6.5,0) {$\sigma$};
\draw  (7.8,.25) to (9.2,.25) to (9.2,-.25) to (7.8,-.25) to (7.8,.25);
\node at (8.5,0) {$\sigma$};
\draw [middleprightarrow] (6,-.25) to (6,-1);
\draw [rightpleftarrow] (7,-.25) to (8.5,-1);
\node at (7.5,-1) {...};
\draw [rightprightarrow] (8,-.25) to (6.5,-1);
\draw [middlepleftarrow] (9,-.25) to (9,-1);
\end{tikzpicture}
\eeq
We denote by $alt_m^{(2)} = (m!)^{-1}\sum_{\sigma\in S_m} sign(\sigma)\sigma^{(2)}$ the projector on $\Lambda^{m}(V\ot V^*)$.
\\

The following is a classical result \cite{ko} (see also \cite{it}).
\bpr
The algebra of invariants 
$$\Lambda^*(\gl(V))^{\gl(V)}\ \simeq\ \Lambda(x_1,x_3,...,x_{2d-1})$$
is the exterior algebra with generators
\beq\lb{giv}
x_{m} = alt^{(2)}_{m}\left(\ 
\begin{tikzpicture}
\draw [middlepleftarrow] (0,0) to [out=90,in=180] (1.5,1) to [out=0,in=90] (3,0);
\draw [middleprightarrow] (.38,0) to [out=90,in=180] (.6,.2) to [out=0,in=90] (.82,0);
\node at (1.5,0) {...};
\draw [middleprightarrow] (2.18,0) to [out=90,in=180] (2.4,.2) to [out=0,in=90] (2.62,0);
\end{tikzpicture}
\ \right) \ \in \ \Lambda^{m}(V\ot V^*)\ .
\eeq
\epr

\bre\lb{ker}
The elements $x_{m}$ defined in \eqref{giv} make sense for all values of $m$.
It can be verified (see \cite{it}) that $x_{m}=0$ for even $m$ or for $m>2d-1$. 
\ere

The natural $\gl(V)$-action on $V$, written as $V\ot V^*\ot V\to V$, takes the shape
$$
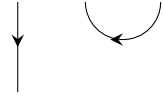
\begin{tikzpicture}
\draw [middleprightarrow] (-0.4,1) to (-0.4,-0.2);
\draw [middlepleftarrow] (0.5,1) to [out=-90,in=180] (1,0.5) to [out=0,in=-90] (1.5,1);
\end{tikzpicture}
$$
The $\gl(V)^{\ot m}$-action on $V^{\ot m}$, written as $(V\ot V^*)^{\ot m}\ot V^{\ot m}\to V^{\ot m}$, has the graphical form
$$
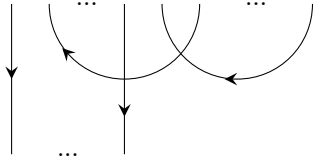
\begin{tikzpicture}
\draw [middleprightarrow] (0,1) to (0,-1);
\draw [middlepleftarrow] (0.5,1) to [out=-90,in=180] (1.5,0);
\draw (1.5,0) to [out=0,in=-90] (2.5,1);
\node at (1,1) {...};
\draw  (1.5,1) to (1.5,0);
\draw [middleprightarrow] (1.5,0) to (1.5,-1);
\node at (3.25,1) {...};
\draw [middlepleftarrow] (2,1) to [out=-90,in=180] (3,0) to [out=0,in=-90] (4,1);
\node at (.75,-1) {...};
\end{tikzpicture}
$$
In particular, the graphical presentation of the action of $x_{m}$ is
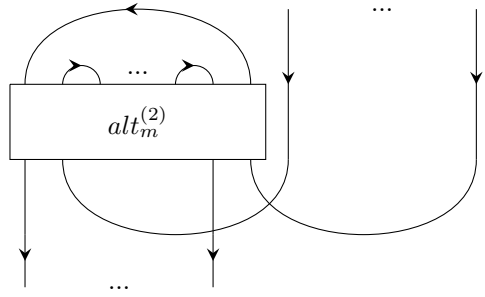
\beq\lb{acti}
\begin{tikzpicture}
\draw [middleprightarrow] (-.5,2) to [out=90,in=180] (-.25,2.25) to [out=0,in=90] (0,2);
\node at (.5,2.15) {...};
\draw [middleprightarrow] (1,2) to [out=90,in=180] (1.25,2.25) to [out=0,in=90] (1.5,2);
\draw [middlepleftarrow] (-1,2) to [out=90,in=180]  (.5,3) to [out=0,in=90] (2,2);
\draw  (-1.2,2) to (2.2,2) to (2.2,1) to (-1.2,1) to (-1.2,2);
\node at (.5,1.5) {$alt_m^{(2)}$};
\draw  (-1,1) to (-1,0);
\draw [middleprightarrow] (-1,0) to (-1,-.7);
\draw  (-.5,1) to [out=-90,in=180]  (1,0) to [out=0,in=-90] (2.5,1);
\draw  (1.5,1) to (1.5,0);
\draw [middleprightarrow] (1.5,0) to (1.5,-.7);
\draw (2,1) to [out=-90,in=180] (3.5,0) to [out=0,in=-90] (5,1);
\node at (.25,-.7) {...};
\draw [middleprightarrow] (2.5,3) to (2.5,1);
\draw [middleprightarrow] (5,3) to (5,1);
\node at (3.75,3) {...};
\end{tikzpicture}
\eeq
The next proposition shows that the $x_{m}$-action on $V^{\ot m}$ coincides up to a scalar multiple with the effect of the element $e_{m}\in k[S_{m}]$ from example \ref{exa}.
\bpr\lb{kox}
The action of $x_{m}$ on $V^{\ot m}$ coincides with $((m-1)!)^{-1}e_{m}$. 
\epr
\bpf
The $x_{m}$-action \eqref{acti} is the alternating sum of terms 
$$
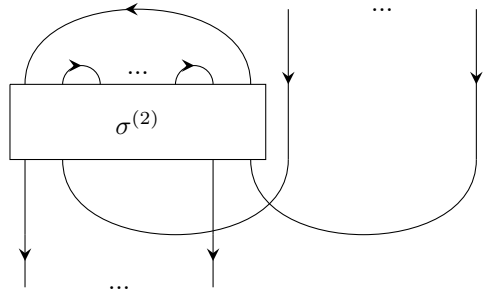
\begin{tikzpicture}
\draw [middleprightarrow] (-.5,2) to [out=90,in=180] (-.25,2.25) to [out=0,in=90] (0,2);
\node at (.5,2.15) {...};
\draw [middleprightarrow] (1,2) to [out=90,in=180] (1.25,2.25) to [out=0,in=90] (1.5,2);
\draw [middlepleftarrow] (-1,2) to [out=90,in=180]  (.5,3) to [out=0,in=90] (2,2);
\draw  (-1.2,2) to (2.2,2) to (2.2,1) to (-1.2,1) to (-1.2,2);
\node at (.5,1.5) {$\sigma^{(2)}$};
\draw  (-1,1) to (-1,0);
\draw [middleprightarrow] (-1,0) to (-1,-.7);
\draw  (-.5,1) to [out=-90,in=180]  (1,0) to [out=0,in=-90] (2.5,1);
\draw  (1.5,1) to (1.5,0);
\draw [middleprightarrow] (1.5,0) to (1.5,-.7);
\draw (2,1) to [out=-90,in=180] (3.5,0) to [out=0,in=-90] (5,1);
\node at (.25,-.7) {...};
\draw [middleprightarrow] (2.5,3) to (2.5,1);
\draw [middleprightarrow] (5,3) to (5,1);
\node at (3.75,3) {...};
\end{tikzpicture}
$$
for $\sigma\in S_m$, which according to \eqref{dup} can be rewritten as 
$$
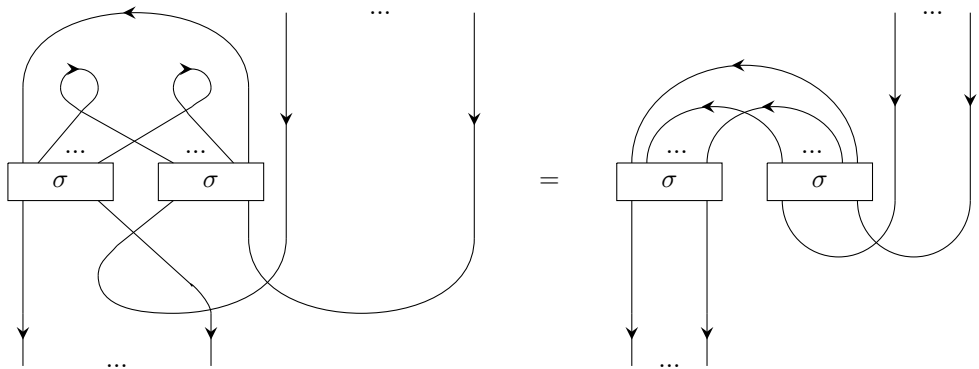
\begin{tikzpicture}
\node at (-.3,2.15) {...};
\node at (1.3,2.15) {...};
\draw [middlepleftarrow] (-1,3) to [out=90,in=180]  (.5,4) to [out=0,in=90] (2,3);
\draw  (-1,3) to (-1,2);
\draw  (2,3) to (2,2);
\draw [middleprightarrow] (-.5,3) to [out=90,in=180] (-.25,3.25) to [out=0,in=90] (0,3);
\draw  (0,3) to [out=-90,in=55] (-.25,2.65);
\draw  (-.25,2.65) to (-.8,2);
\draw  (-.5,3) to [out=-90,in=140] (-.25,2.7);
\draw  (-.25,2.7) to (1,2);
\draw [middleprightarrow] (1,3) to [out=90,in=180] (1.25,3.25) to [out=0,in=90] (1.5,3);
\draw  (1.5,3) to [out=-90,in=45] (1.25,2.75);
\draw  (1.25,2.75) to (0,2);
\draw  (1,3) to [out=-90,in=120] (1.2,2.65);
\draw  (1.2,2.65) to (1.8,2);
\draw  (-1.2,2) to (0.2,2) to (0.2,1.5) to (-1.2,1.5) to (-1.2,2);
\node at (-.5,1.75) {$\sigma$};
\draw  (.8,2) to (2.2,2) to (2.2,1.5) to (.8,1.5) to (.8,2);
\node at (1.5,1.75) {$\sigma$};
\draw  (-1,1.5) to (-1,0);
\draw [middleprightarrow] (-1,0) to (-1,-.7);
\draw  (0,.5) to [out=-90,in=180]  (1,0) to [out=0,in=-90] (2.5,1);
\draw  (0,.5) to [out=90,in=-130]  (.25,.9);
\draw  (.25,.9) to (1,1.5);
\draw  (0,1.5) to (1.25,.37); 
\draw  (1.25, .357) to [out=127,in=90] (1.5,0); 
\draw [middleprightarrow] (1.5,0) to (1.5,-.7);
\draw (2,1) to (2,1.5);
\draw (2,1) to [out=-90,in=180] (3.5,0) to [out=0,in=-90] (5,1);
\node at (.25,-.7) {...};
\draw [middleprightarrow] (2.5,4) to (2.5,1);
\draw [middleprightarrow] (5,4) to (5,1);
\node at (3.75,4) {...};
\node at (6,1.75) {=};
\end{tikzpicture}
\hspace{1.5pc}
\begin{tikzpicture}
\node at (-.4,2.15) {...};
\node at (1.4,2.15) {...};
\draw [middlepleftarrow] (-1,2) to [out=90,in=180]  (.5,3.3) to [out=0,in=90] (2,2);
\draw [middlepleftarrow] (-.8,2) to [out=90,in=180] (.1,2.75) to [out=0,in=90] (1,2);
\draw [middlepleftarrow] (0,2) to [out=90,in=180] (.9,2.75) to [out=0,in=90] (1.8,2);
\draw  (-1.2,2) to (0.2,2) to (0.2,1.5) to (-1.2,1.5) to (-1.2,2);
\node at (-.5,1.75) {$\sigma$};
\draw  (.8,2) to (2.2,2) to (2.2,1.5) to (.8,1.5) to (.8,2);
\node at (1.5,1.75) {$\sigma$};
\draw  (-1,1.5) to (-1,0);
\draw [middleprightarrow] (-1,0) to (-1,-.7);
\node at (-.5,-.7) {...};
\draw  (0,1.5) to (0,0);
\draw [middleprightarrow] (0,0) to (0,-.7);
\draw (1,1.5) to [out=-90,in=180] (1.75,0.75) to [out=0,in=-90] (2.5,1.5);
\draw (2,1.5) to [out=-90,in=180] (2.75,0.75) to [out=0,in=-90] (3.5,1.5);
\draw [middleprightarrow] (2.5,4) to (2.5,1.5);
\draw [middleprightarrow] (3.5,4) to (3.5,1.5);
\node at (3,4) {...};
\end{tikzpicture}
$$
The identity \eqref{inp} allows us to rewrite it further as follows
$$
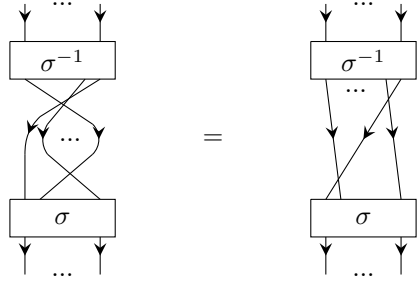
\begin{tikzpicture}
\draw [middleprightarrow] (-1,4.6) to (-1,4.1);
\node at (-.5,4.6) {...};
\draw [middleprightarrow] (0,4.6) to (0,4.1);
\draw  (-1.2,4.1) to (0.2,4.1) to (0.2,3.6) to (-1.2,3.6) to (-1.2,4.1);
\node at (-.5,3.85) {$\sigma^{-1}$};
\draw  (-.8,3.1) to (0,3.6);
\draw [middleprightarrow] (-.8,3.1) to [out=-135,in=90] (-1,2.6);
\draw  (-1,2.6) to (-1,2);
\draw  (-.1,3) to (-1,3.6); 
\draw [middleprightarrow] (-.1,3) to [out=-45,in=45] (-.1,2.6);
\draw  (-.1,2.6) to (-.8,2); 
\draw  (-.7,3) to (-.2,3.6); %
\draw [middleprightarrow] (-.7,3) to [out=-120,in=125] (-.7,2.6);
\draw  (-.7,2.6) to (0,2); %
\draw  (-1.2,2) to (0.2,2) to (0.2,1.5) to (-1.2,1.5) to (-1.2,2);
\node at (-.5,1.75) {$\sigma$};
\draw [middleprightarrow] (-1,1.5) to (-1,1);
\node at (-.5,1) {...};
\draw [middleprightarrow] (0,1.5) to (0,1);
\node at (-.4,2.8) {...};
\node at (1.5,2.8) {=};
\draw [middleprightarrow] (3,4.6) to (3,4.1);
\node at (3.5,4.6) {...};
\draw [middleprightarrow] (4,4.6) to (4,4.1);
\draw  (2.8,4.1) to (4.2,4.1) to (4.2,3.6) to (2.8,3.6) to (2.8,4.1);
\node at (3.5,3.85) {$\sigma^{-1}$};
\node at (3.4,3.45) {...};
\draw  [middleprightarrow] (4,3.6) to (3,2);
\draw  [middleprightarrow] (3,3.6) to (3.2,2);
\draw  [middleprightarrow] (3.8,3.6) to (4,2);
\draw  (2.8,2) to (4.2,2) to (4.2,1.5) to (2.8,1.5) to (2.8,2);
\node at (3.5,1.75) {$\sigma$};
\draw [middleprightarrow] (3,1.5) to (3,1);
\node at (3.5,1) {...};
\draw [middleprightarrow] (4,1.5) to (4,1);
\end{tikzpicture}
$$
Thus the $x_{m}$-action \eqref{acti} is the alternating sum of the conjugates of the long cycle in $S_m$.
Since the centraliser of the long cycle in $S_m$ has order $m$ we get the answer $((m-1)!)^{-1}e_{m}$. 
\epf

The assignment $[1]\mapsto V$ extends to a symmetric tensor functor
$$SW:\S\to \Rep(\gl(V))\ ,$$
the {\em Schur-Weyl functor} (see \cite{dm} for details). 
According to \eqref{hmc} we have a pair of homomorphisms
$$H^*(\S)\ \to\ H^*(SW)\ \leftarrow\ H^*(\Rep(\gl(V)))\ .$$
The following theorem shows that the first homomorphism is surjective and the second homomorphism is an isomorphism. 
\bth\lb{cswf}
The deformation cohomology of the Schur-Weyl functor $SW$ is the exterior algebra 
$$H^*(SW) = \Lambda(e_1,e_3,...,e_{2d-1}) .$$ 
The homomorphism $H^*(\S)\ \to\ H^*(SW)$ is the quotioning by the ideal generated by $e_{s}, s>2d-1$.
\nl
The homomorphism $H^*(\Rep(\gl(V)))\ \to\ H^*(SW)$ is an isomorphism sending $x_m$ to $((m-1)!)^{-1}e_{m}$. 
\eth
\bpf
The classical Schur-Weyl duality implies that the Schur-Weyl functor is full.
Moreover it identifies the kernel of (i.e. the ideal of morphisms in $\S$ annihilated by) the Schur-Weyl functor $SW$ with the tensor ideal $J$ generated by $alt_{d+1}$. By lemma \ref{cef} the deformation cohomology $H^*(SW)$ of the Schur-Weyl functor $SW$ coincides with the  deformation cohomology $H^*(\S/J)$ of the quotient category $\S/J$ of $\S$ by this ideal.
\nl
The quotient $\S/J$ is  a Schur-Weyl category $\SW(A_*)$ of the multiplicative sequence $A_n = k[S_n]/J_n$. 
The methods of section \ref{cdi} also apply to this sequence. Namely the cohomology of horisontal complexes is concentrated in the top degree and coincides with the deformation cohomology $H^*(\S/J)$.
In particular the homomorphism $H^*(\S)\to H^*(\S/J)=H^*(SW)$ is surjective.
\nl
Proposition \ref{kox} says that the image of $s_s\in H^*(\S)$ is $H^*(SW)$ coincides up to a scalar with the image of $x_s\in H^*(\Rep(\gl(V)))$. This shows that the homomorphism $H^*(\Rep(\gl(V)))\ \to\ H^*(SW)$ is an isomorphism.
Finally remark \ref{ker} describes the kernel of the homomorphism $H^*(\S)\ \to\ H^*(SW)$. 
\epf

\subsection{Deformation cohomology of $\Rep(\gl(V)\ot A)$}

Let $A$ be a commutative algebra and $\g$ be a Lie algebra.
The tensor product $\g\ot A$ has s structure of Lie algebra
$$[x\ot a,y\ot b] = [x,y]\ot ab,\qquad x,y\in\g,\ a,b\in A\ .$$

The following is an auxiliary result.
For $a\in A$ we will denote by $a^{(i)} = 1\ot...\ot 1\ot a\ot 1...\ot 1$ the decomposable tensor in $A^{\ot n}$ with all by $i$-th components being one.
\ble\lb{ale}
Let $A$ be a commutative algebra with $dim(A)>1$ and such that $A^{\ot n}$ has no zero-divisors. 
Let $M$ be a free $A^{\ot n}$-module.
Let $m_i\in M$ be such that $\sum_{i=1}^n m_i a^{(i)} = 0$ for any $a\in A$.
Then $m_i=0$ for any $i=1,...,n$.
\ele
\bpf
It is enough to prove the lemma for $M=A^{\ot n}$.
Take $a\in A$ to be non-scalar.
The system of equations $\sum_{i=1}^n m_i (a^s)^{(i)} = 0$ for $s=1,...,n$ can be considered as a linear system with unknowns $m_i$.
By multiplying with the adjugate matrix we can see that $m_i$ are annihilated by the determinant of the system.
This determinant is the Vandermonde determinant $\prod_{i<j}(a^{(i)}-a^{(j)})$. 
It is non-zero by the choice of $a$ and we get the conclusion of the lemma by the absence of zero-divisors in $A^{\ot n}$. 
\epf

\bpr
Let $A$ be a commutative algebra with $dim(A)>1$ and such that $A^{\ot n}$ has no zero-divisors for any $n$. 
Let $\g$ be a Lie algebra and let $\z(\g)$ be its centre.
Then the algebra of exterior invariants of $\g\ot A$ is isomorphic to the exterior algebra of $z(\g)\ot A$
$$\Lambda^*(\g\ot A)^{\g\ot A}\ \simeq\ \Lambda^*(z(\g)\ot A)\ .$$
\epr
\bpf
Using the isomorphism $(\g\ot A)^{\ot n}\simeq \g^{\ot n}\ot A^{\ot n}$ we can identify $\Lambda^*(\g\ot A)$ with the subspace of $S_n$-anti-invariants in $\g^{\ot n}\ot A^{\ot n}$. 
An element $u\in \g^{\ot n}\ot A^{\ot n}$ is $\g\ot A$-invariant if 
$$\sum_{i=1}^n ad_x^{(i)}(u) a^{(i)} = 0$$ 
for any $x\in \g$ and $a\in A$.
Here $ad_x(y) = [x,y]$ for $y\in \g$ and $ad_x^{(i)}(u)$ is the result of applying $ad_x$ to the $i$-th component (in the $\g^{\ot n}$ part) of $u$.
The product $ad_x^{(i)}(u) a^{(i)}$ is in the sense of the natural $A^{\ot n}$-module structure on $\g^{\ot n}\ot A^{\ot n}$. 
Now the proposition follows from lemma \ref{ale}. 
\epf

\bre\lb{ina}
By the above proposition (assuming that $A$ is a commutative algebra with $dim(A)>1$ and without zero-divisors)
we have that the algebra of invariants of $\gl(V)\ot A$ is isomorphic to the exterior algebra of $A$
$$\Lambda^*(\gl(V)\ot A)^{\gl(V)\ot A}\ \simeq\ \Lambda^*(A)\ .$$
\ere

The natural $\gl(V)$-module $V$ give rise to a $\gl(V)\ot A$-module structure on $V\ot A$
$$(x\ot a).(v\ot b) = x(v)\ot ab,\qquad x\in \gl(V), v\in V, a,b\in A\ .$$

As before the assignment $[1]\mapsto V$ extends to a symmetric tensor functor
\beq\lb{swfa}SW:\S(A)\to \Rep(\gl(V\ot A))\ ,\eeq
the {\em Schur-Weyl functor} (see \cite{dm} for details). 
Again according to \eqref{hmc} we have a pair of homomorphisms
\beq\lb{swh}H^*(\S(A))\ \to\ H^*(SW)\ \leftarrow\ H^*(\Rep(\gl(V\ot A)))\ .\eeq
The following theorem shows that the first homomorphism is surjective and the second homomorphism is an isomorphism. 
\bth\lb{swda}
Let $A$ be a commutative algebra with $dim(A)>1$ and such that $A^{\ot n}$ has no zero-divisors for any $n$. 
Then the deformation cohomology of the Schur-Weyl functor \eqref{swfa} is the exterior algebra 
$$H^*(SW) = \Lambda^*(A)$$ 
and the homomorphisms \eqref{swh} are isomorphisms. 
\eth
\bpf
The proof is similar to the proof of theorem \ref{cswf}.
The Schur-Weyl duality implies that the Schur-Weyl functor is full.
By lemma \ref{cef} the deformation cohomology $H^*(SW)$ of the Schur-Weyl functor  \eqref{swfa} coincides with the deformation cohomology $H^*(\S(A)/J)$ of the quotient category $\S(A)/J$ of $\S$ by the kernel of the Schur-Weyl functor.
\nl
The quotient $\S(A)/J$ is a Schur-Weyl category $\SW(A_*)$ of the multiplicative sequence $A_n = (A^{\ot n}*S_n)/J_n$. 
The methods of section \ref{cdi} also apply to this sequence. Namely the cohomology of the horisontal complexes is concentrated in the top degree and coincides with the deformation cohomology $H^*(\S(A)/J)$.
In particular the homomorphism $H^*(\S(A))\to H^*(\S(A)/J)=H^*(SW)$ is surjective.
\nl
By theorem \ref{cfs} the cohomology $H^*(\S(A))$ is generated by degree one elements and so is $H^*(SW)$.
Now by remark \ref{ina} the cohomology $H^*(\Rep(\gl(V\ot A)))$ is also generated by degree one elements and the space of generators is the same as for $\S(A)$.
Comparing (with the help of remark \ref{pfs}) the effect of generators on $V\ot A\in \Rep(\gl(V\ot A))$ we get the theorem.
\epf

\end{document}